\documentclass[12pt]{article}
\usepackage{amsmath,amsfonts,amssymb,amsthm}
\input amssym.def
\begin{document}
\def \Z{\Bbb Z}
\def \C{\Bbb C}
\def \R{\Bbb R}
\def \Q{\Bbb Q}
\def \N{\Bbb N}

\def \A{{\mathcal{A}}}
\def \D{{\mathcal{D}}}
\def \E{{\mathcal{E}}}
\def \L{\mathcal{L}}
\def \S{{\mathcal{S}}}
\def \wt{{\rm wt}}
\def \tr{{\rm tr}}
\def \span{{\rm span}}
\def \Res{{\rm Res}}
\def \Der{{\rm Der}}
\def \End{{\rm End}}
\def \Ind {{\rm Ind}}
\def \Irr {{\rm Irr}}
\def \Aut{{\rm Aut}}
\def \GL{{\rm GL}}
\def \Hom{{\rm Hom}}
\def \mod{{\rm mod}}
\def \ann{{\rm Ann}}
\def \ad{{\rm ad}}
\def \rank{{\rm rank}\;}
\def \<{\langle}
\def \>{\rangle}

\def \g{{\frak{g}}}
\def \h{{\hbar}}
\def \k{{\frak{k}}}
\def \sl{{\frak{sl}}}
\def \gl{{\frak{gl}}}

\def \be{\begin{equation}\label}
\def \ee{\end{equation}}
\def \bex{\begin{example}\label}
\def \eex{\end{example}}
\def \bl{\begin{lem}\label}
\def \el{\end{lem}}
\def \bt{\begin{thm}\label}
\def \et{\end{thm}}
\def \bp{\begin{prop}\label}
\def \ep{\end{prop}}
\def \br{\begin{rem}\label}
\def \er{\end{rem}}
\def \bc{\begin{coro}\label}
\def \ec{\end{coro}}
\def \bd{\begin{de}\label}
\def \ed{\end{de}}

\newcommand{\m}{\bf m}
\newcommand{\n}{\bf n}
\newcommand{\nno}{\nonumber}
\newcommand{\nord}{\mbox{\scriptsize ${\circ\atop\circ}$}}
\newtheorem{thm}{Theorem}[section]
\newtheorem{prop}[thm]{Proposition}
\newtheorem{coro}[thm]{Corollary}
\newtheorem{conj}[thm]{Conjecture}
\newtheorem{example}[thm]{Example}
\newtheorem{lem}[thm]{Lemma}
\newtheorem{rem}[thm]{Remark}
\newtheorem{de}[thm]{Definition}
\newtheorem{hy}[thm]{Hypothesis}
\makeatletter
\@addtoreset{equation}{section}
\def\theequation{\thesection.\arabic{equation}}
\makeatother
\makeatletter

\begin{center}
{\Large \bf Some quantum vertex algebras of Zamolodchikov-Faddeev
type}
\end{center}

\begin{center}
{Martin Karel and Haisheng Li\footnote{Partially supported
by NSA grant H98230-05-1-0018 and NSF grant DMS-0600189}\\
Department of Mathematical Sciences, Rutgers University, Camden, NJ
08102}
\end{center}

\begin{abstract}
This is a continuation of a previous study of quantum vertex
algebras of Zamolodchikov-Faddeev type. In this paper, we focus our
attention on the special case associated to diagonal unitary
rational quantum Yang-Baxter operators. We prove that the associated
weak quantum vertex algebras, if not zero, are irreducible quantum
vertex algebras with a normal basis in a certain sense.
\end{abstract}

\section{Introduction}
In an attempt to associate vertex algebra-like structures with
quantum affine algebras, in \cite{li-qva1} we studied general formal
vertex operators (=quantum fields) on an {\em arbitrary} vector
space, to see what kind of algebraic structures they could possibly
``generate''. More specifically, we studied what we called ``quasi
compatible'' sets of formal vertex operators on a vector space and
we proved that any quasi compatible set generates a nonlocal vertex
algebra in a certain canonical way. (Whereas vertex algebras are
analogues of commutative and associative algebras, nonlocal vertex
algebras are analogues of noncommutative associative algebras.) It
follows from this general result that to most of the interesting
algebras one can associate nonlocal vertex algebras. For example,
for any highest weight module $W$ for a quantum affine algebra, the
generating functions in Drinfeld's realization form a quasi
compatible set of vertex operators on $W$ and therefore they
generate a nonlocal vertex algebra.

As general nonlocal vertex algebras are too ``wild'',  one hopes to
get ``better'' families of nonlocal vertex algebras from ``better''
families of vertex operators. Motivated by the $\S$-locality axiom
in Etingof-Kazhdan's notion of quantum vertex operator algebra
\cite{ek},  we then continued to study what we called $\S$-local
sets and quasi $\S$-local sets of vertex operators. We proved that
from a quasi $\S$-local set of vertex operators, indeed one can get
a better nonlocal vertex algebra whose adjoint vertex operators form
an $\S$-local set of vertex operators (on the nonlocal vertex
algebra). Motivated by this result and by Etingof-Kazhdan's notion
of quantum vertex operator algebra we then formulated and studied a
notion of weak quantum vertex algebra and a notion of quantum vertex
algebra.

 A {\em weak quantum vertex
algebra} (over $\C$) is a vector space $V$ equipped with a
distinguished vector ${\bf 1}$ and a linear map $Y: V\rightarrow
\Hom(V,V((x)))$ satisfying the condition that $Y({\bf 1},x)=1$ (the
identity operator on $V$),
$$Y(v,x){\bf 1}\in V[[x]]\;\;\mbox{ and }\ \ (Y(v,x){\bf
1})|_{x=0}=v\ \ \ \mbox{ for }v\in V,$$ and that
 for any $u,v\in V$, there exist
$u^{(i)},v^{(i)}\in V,\; f_{i}(x)\in \C((x)),\; i=1,\dots,r$, such
that
\begin{eqnarray}\label{e1.1}
& &x_{0}^{-1}\delta\left(\frac{x_{1}-x_{2}}{x_{0}}\right)
Y(u,x_{1})Y(v,x_{2})\nonumber\\
& &\ \ \ \ \
-x_{0}^{-1}\delta\left(\frac{x_{2}-x_{1}}{-x_{0}}\right)\sum_{i=1}^{r}
f_{i}(-x_{0})Y(v^{(i)},x_{2})Y(u^{(i)},x_{1})\nonumber\\
&=&x_{2}^{-1}\delta\left(\frac{x_{1}-x_{0}}{x_{2}}\right)
Y(Y(u,x_{0})v,x_{2}).
\end{eqnarray}
A {\em quantum vertex algebra} is a weak quantum vertex algebra $V$
equipped with a unitary rational quantum Yang-Baxter operator (with
one parameter)
$$\S(x): V\otimes
V\rightarrow V\otimes V\otimes \C((x))$$ such that for any $u,v\in
V$, (\ref{e1.1}) holds with
$$\S(x)(v\otimes u)=\sum_{i=1}^{r}v^{(i)}\otimes u^{(i)}\otimes
f_{i}(x)$$ and such that some other conditions hold. The notion of
quantum vertex algebra generalizes the notions of vertex algebra and
vertex superalgebra.

The designation of the notion of weak quantum vertex algebra was
motivated by the following conceptual result. Let $W$ be any vector
space and set $\E(W)=\Hom (W,W((x)))$.
 A subset $T$ of $\E(W)$ is said
to be {\em $\S$-local} if for any $a(x),b(x)\in T$, there exist
(finitely many)
$$a^{(i)}(x),b^{(i)}(x)\in T,\; f_{i}(x)\in \C((x)),\; i=1,\dots,r$$
such that
\begin{eqnarray}
(x_{1}-x_{2})^{k}a(x_{1})b(x_{2}) =
(x_{1}-x_{2})^{k}\sum_{i=1}^{r}f_{i}(x_{2}-x_{1})b^{(i)}(x_{2})a^{(i)}(x_{1})
\end{eqnarray}
for some nonnegative integer $k$. Let $T$ be an $\S$-local subset of
$\E(W)$. For $a(x),b(x)\in T$ with the above information, we define
\begin{eqnarray}
& &Y_{\E}(a(x),x_{0})b(x)
=\Res_{x_{1}}x_{0}^{-1}\delta\left(\frac{x_{1}-x}{x_{0}}\right)
a(x_{1})b(x)\nonumber\\
& &\hspace{2cm} \ \ \
-\Res_{x_{1}}x_{0}^{-1}\delta\left(\frac{x-x_{1}}{-x_{0}}\right)
\sum_{i=1}^{r}f_{i}(x_{2}-x_{1})b^{(i)}(x_{2})a^{(i)}(x_{1}).
\end{eqnarray}
It was proved in \cite{li-qva1} that any $\S$-local subset of
$\E(W)$ generates a weak quantum vertex algebra with $W$ as a
canonical module. This particular result generalizes a result of
\cite{li-local}, which states that for any vector space $W$, every
local subset of $\E(W)$ generates a vertex algebra with $W$ as a
module (see \cite{li-twisted}, \cite{li-g1}, \cite{li-gamma}, and
\cite{li-infinity} for similar results).

Regarding quantum vertex algebras, a variant, which was formulated
in \cite{li-qva1}, of (\cite{ek}, Proposition 1.11),  is that if a
weak quantum vertex algebra $V$ is nondegenerate in the sense of
\cite{ek}, $V$ is a quantum vertex algebra with $\S(x)$ uniquely
determined. Furthermore, it was proved in \cite{li-qva2} (cf.
\cite{li-simple} for vertex algebras) that if a nonlocal vertex
algebra $V$ (over $\C$) is of countable dimension and if $V$ as a
(left) $V$-module is irreducible, then $V$ is nondegenerate.

In \cite{li-qva2}, to demonstrate how to use the general
construction theorem which was established in \cite{li-qva1}, we
constructed some quantum vertex algebras by using algebras of
Zamolodchikov-Faddeev type (see \cite{za}, \cite{fa}). More
specifically, let $H$ be a finite-dimensional vector space equipped
with a bilinear form $\<\cdot,\cdot\>$ and let $\S(x): H\otimes
H\rightarrow H\otimes H\otimes \C[[x]]$ be a linear map. We studied
weak quantum vertex algebras $V$ equipped with a linear map $\phi:
H\rightarrow V$ such that $V$ is generated by $\phi(H)$, satisfying
the relations for $u,v\in H$:
\begin{eqnarray*}
& &Y(\phi(u),x_{1})Y(\phi(v),x_{2})-\sum_{i=1}^{r}f_{i}(x_{2}-x_{1})
Y(\phi(v^{(i)}),x_{2})Y(\phi(u^{(i)}),x_{1})\\
&=&\<u,v\>x_{2}^{-1}\delta\left(\frac{x_{1}}{x_{2}}\right),
\end{eqnarray*}
where $\S(x)(v\otimes u)=\sum_{i=1}^{r}v^{(i)}\otimes u^{(i)}\otimes
f_{i}(x)$.  We then constructed a universal weak quantum vertex
algebra $V(H,\S)$ and proved that $V(H,\S)$ is nondegenerate for
some special cases, so that $V(H,\S)$ is a quantum vertex algebra
(see Section 4).  (If $\phi$ is injective and if $V(H,\S)$ is
nondegenerate, then $\S(x)$ is necessarily a unitary rational
quantum Yang-Baxter operator.)

This current paper is a natural continuation of our studies on weak
quantum vertex algebras $V(H,\S)$. In this paper, we focus our
attention on the special case in which $\S(x)$ is ``diagonal''. Let
${\bf Q}(x)=(q_{ij}(x))$ be an $l\times l$ matrix in $\C[[x]]$ with
$q_{ij}(x)q_{ji}(-x)=1$ for $1\le i,j\le l$. Take $H$ to be a
$2l$-dimensional vector space with a basis $\{ u^{(i)},v^{(i)}\;|\;
1\le i\le l\}$, equip $H$ with the bilinear form $\<\cdot,\cdot\>$
defined by
$$\<u^{(i)},u^{(j)}\>=0=\<v^{(i)},v^{(j)}\>,\ \ \
\<u^{(i)},v^{(j)}\>=\delta_{ij}=-q_{ii}(0)\<v^{(j)},u^{(i)}\>   $$
for $1\le i,j\le l$,  and define $\S(x): H\otimes H\rightarrow
H\otimes H\otimes \C[[x]]$ by
\begin{eqnarray*}
& &\S(x)(u^{(j)}\otimes u^{(i)})=u^{(j)}\otimes u^{(i)}\otimes
q_{ij}(x),\ \ \ \S(x)(v^{(j)}\otimes v^{(i)})=v^{(j)}\otimes
v^{(i)}\otimes
q_{ij}(x),\\
& &\S(x)(v^{(j)}\otimes u^{(i)})=v^{(j)}\otimes u^{(i)}\otimes
q_{ji}(-x),\ \ \S(x)(u^{(j)}\otimes v^{(i)})=u^{(j)}\otimes
v^{(i)}\otimes q_{ji}(-x)
\end{eqnarray*}
for $1\le i,j\le l$. For this pair $(H,\S)$, the associated
universal weak quantum vertex algebra $V(H,\S)$  is alternatively
denoted by $V_{{\bf Q}(x)}$.  Our main goal is to determine the
structure and establish the nondegeneracy of $V_{{\bf Q}(x)}$. As
one of the main results in this paper, we prove that $V_{{\bf
Q}(x)}$ is either zero or irreducible as a $V_{{\bf Q}(x)}$-module
with a normal basis in a certain sense, so that it is a
nondegenerate quantum vertex algebra.

As a technical step, we first study the case with ${\bf Q}(x)={\bf
Q}=(q_{ij})$ a complex matrix with $q_{ij}q_{ji}=1$ for $1\le i,j\le
l$. We define an associative algebra $\A_{\bf Q}$ with generators
$X_{i,m},Y_{i,m}$ for $1\le i\le l,\; m\in \Z$, subject to relations
\begin{eqnarray*}
& &X_{i,m}X_{j,n}=q_{ij}X_{j,n}X_{i,m}, \ \ \
Y_{i,m}Y_{j,n}=q_{ij}Y_{j,n}Y_{i,m},\\
& &\ \ \ \
X_{i,m}Y_{j,n}-q_{ji}Y_{j,n}X_{i,m}=\delta_{ij}\delta_{m+n+1,0}
\end{eqnarray*}
for $1\le i,j\le l,\; m,n\in \Z$.  If $q_{ij}=1$ for all
$i,j=1,\dots,l$, $\A_{\bf Q}$ is a Weyl algebra and if $q_{ij}=-1$
for all $i,j=1,\dots,l$, $\A_{\bf Q}$ is a Clifford algebra. (In
general, we have $q_{ii}=\pm 1$ for $1\le i\le l$, so that $\A_{\bf
Q}$ contains Weyl algebras or Clifford algebras as subalgebras.)
 Let $V_{\bf Q}$ denote the quotient
module of the regular left $\A_{\bf Q}$-module modulo the left ideal
generated by $X_{i,m},Y_{i,m}$ for $1\le i\le l,\; m\ge 0$.  It is
proved that $V_{\bf Q}$ is irreducible and there is a canonical
simple quantum vertex algebra structure on $V_{\bf Q}$. Furthermore,
$V_{\bf Q}$ has a conformal vector of central charge $l$. Note that
 quantum coordinate algebras, closely related to $\A_{\bf Q}$, have
appeared before in the study of noncommutative geometry (cf.
\cite{manin1}, \cite{manin2}) and physics noncommutative field
theory (cf. \cite{ku}). Our emphasis in this study is on quantum
vertex algebras and modules.

To achieve our results for the nonconstant case with ${\bf
Q}(x)=(q_{ij}(x))$, we make use of a certain filtration. We define
an increasing filtration $F=\{ F_{n}\}_{n\ge 0}$ of $V_{{\bf Q}(x)}$
by using the canonical generators and we prove that if the weak
quantum vertex algebra $V_{{\bf Q}(x)}$ is nonzero, the associated
graded weak quantum vertex algebra ${\rm Gr}_{F}(V_{{\bf Q}(x)})$ is
isomorphic to the irreducible quantum vertex algebra $V_{{\bf
Q}(0)}$. Then we prove that $V_{{\bf Q}(x)}$, if not zero, is an
irreducible quantum vertex algebra.
 Furthermore, we show that for certain cases, $V_{{\bf
Q}(x)}$ is indeed nonzero.

One of the authors (H. Li) would like to thank Igor Frenkel for a
discussion on quantum vertex algebras and Zamolodchikov-Faddeev
algebras. Some time ago, he and J. Ding had some similar ideas. Part
of this paper was done during H.Li's visit in January 2007 at the
Chern Institute of Mathematics, Nankai University, Tianjin, China.
We would like to thank Professors Chengming Bai and Weiping Zhang
for their hospitality.

This paper is organized as follows: In Section 2, we present some
basic results about general nonlocal vertex algebras. In Section 3,
we study the associative algebras $\A_{\bf Q}$ and quantum vertex
algebras $V_{\bf Q}$. In Section 4, we study quantum vertex algebras
$V_{{\bf Q}(x)}$.

\section{Some results for general nonlocal vertex algebras}
In this section we recall the notions of nonlocal vertex algebra,
weak quantum vertex algebra, and quantum vertex algebra, and we
present some basic results on increasing filtrations for general
nonlocal vertex algebras associated with a generating subset.

We start by recalling the definition of a nonlocal vertex algebra
{}from \cite{li-qva1} (cf. \cite{kac}, \cite{bk}, \cite{li-g1}). A
{\em nonlocal vertex algebra} is a vector space $V$, equipped with a
linear map
\begin{eqnarray}
Y: V &\rightarrow & \Hom (V,V((x)))\subset (\End V)[[x,x^{-1}]]\nonumber\\
v&\mapsto& Y(v,x)=\sum_{n\in \Z}v_{n}x^{-n-1}\;\;\; (v_{n}\in \End V)
\end{eqnarray}
and equipped with a distinguished vector ${\bf 1}$, satisfying the
condition that for $v\in V$
\begin{eqnarray}
& &Y({\bf 1},x)v=v,\\
& &Y(v,x){\bf 1}\in V[[x]]\;\;\mbox{ and }\;\; \lim_{x\rightarrow
0}Y(v,x){\bf 1}=v
\end{eqnarray}
and the condition that for $u,v,w\in V$, there exists a nonnegative
integer $l$ such that
\begin{eqnarray}
(x_{0}+x_{2})^{l}Y(u,x_{0}+x_{2})Y(v,x_{2})w=
(x_{0}+x_{2})^{l}Y(Y(u,x_{0})v,x_{2})w.
\end{eqnarray}
Following \cite{ek}, let $Y(x): V\otimes V\rightarrow V((x))$ be the
linear map associated with $Y$.

For a nonlocal vertex algebra $V$, let $\D$ be the linear operator
on $V$ defined by  $\D (v)=v_{-2}{\bf 1}$ for $v\in V$. Then
\begin{eqnarray}
[\D, Y(v,x)]=Y(\D v,x)=\frac{d}{dx}Y(v,x)\ \ \ \mbox{ for }v\in V.
\end{eqnarray}
Furthermore, we have $Y(v,x){\bf 1}=e^{x\D}v$  and $\D ({\bf 1})=0$
for $v\in V$.

For a nonlocal vertex algebra $V$, a {\em $V$-module} is a vector
space $W$ equipped with a linear map $Y_{W}: V\rightarrow \Hom
(W,W((x)))$ satisfying the condition that
$$Y_{W}({\bf 1},x)=1_{W}\ \ \mbox{
(the identity operator on } W)$$
 and for $u,v\in V,\; w\in W$, there
exists $l\in \N$ such that
$$(x_{0}+x_{2})^{l}Y_{W}(u,x_{0}+x_{2})Y_{W}(v,x_{2})w=
(x_{0}+x_{2})^{l}Y_{W}(Y(u,x_{0})v,x_{2})w.$$

Next, we recall from \cite{li-qva1} the notions of weak quantum
vertex algebra and quantum vertex algebra. A {\em weak quantum
vertex algebra} is a vector space $V$ (over $\C$) equipped with a
distinguished vector ${\bf 1}$ and a linear map $Y$ from $V$ to
$\Hom (V,V((x))))$, satisfying the condition that
\begin{eqnarray}
& &Y({\bf 1},x)=1,\\
& &Y(v,x){\bf 1}\in V[[x]]\;\;\mbox{ and }\;\;\lim_{x\rightarrow
0}Y(v,x){\bf 1}=v
\end{eqnarray}
for $v\in V$, and that for any $u,v\in V$, there exist
$u^{(i)},v^{(i)}\in V,\; f_{i}(x)\in \C((x)),\; i=1,\dots,r$, such
that
\begin{eqnarray}\label{eS-jacobi}
& &x_{0}^{-1}\delta\left(\frac{x_{1}-x_{2}}{x_{0}}\right)Y(u,x_{1})Y(v,x_{2})
\nonumber\\
& &\ \ \ -x_{0}^{-1}\delta\left(\frac{x_{2}-x_{1}}{-x_{0}}\right)\sum_{i=1}^{r}
f_{i}(-x_{0})Y(v^{(i)},x_{2})Y(u^{(i)},x_{1})\nonumber\\
&=&x_{2}^{-1}\delta\left(\frac{x_{1}-x_{0}}{x_{2}}\right)
Y(Y(u,x_{0})v,x_{2}).
\end{eqnarray}
Alternatively, a weak quantum vertex algebra is a nonlocal vertex
algebra that satisfies $\S$-locality (cf. \cite{ek}) in the sense
that for any $u,v\in V$, there exist $u^{(i)},v^{(i)}\in V,\;
f_{i}(x)\in \C((x)),\; i=1,\dots,r$, such that
\begin{eqnarray}
(x_{1}-x_{2})^{k}Y(u,x_{1})Y(v,x_{2})
=(x_{1}-x_{2})^{k}\sum_{i=1}^{r}f_{i}(x_{2}-x_{1})
Y(v^{(i)},x_{2})Y(u^{(i)},x_{1})
\end{eqnarray}
for some nonnegative integer $k$ (depending on $u$ and $v$).

The notion of quantum vertex algebra involves a quantum Yang-Baxter
operator. A {\em unitary rational quantum Yang-Baxter operator} with
one parameter on a vector space $U$ is a linear map
$$\S(x): U\otimes U\rightarrow U\otimes U\otimes \C((x))$$ such that
\begin{eqnarray*}
\S(x)\S^{21}(-x)&=&1,\\
\S^{12}(x_{1})\S^{13}(x_{1}-x_{2})\S^{23}(x_{2})
&=&\S^{23}(x_{2})\S^{13}(x_{1}-x_{2})\S^{12}(x_{1}).
\end{eqnarray*}
A {\em quantum vertex algebra} (cf. \cite{ek}) is a weak quantum
vertex algebra $V$ equipped with a unitary rational quantum
Yang-Baxter operator
$$\S(x): V\otimes V\rightarrow V\otimes V\otimes \C((x))$$
such that for $u,v\in V$, (\ref{eS-jacobi})  holds with
$\S(x)(v\otimes u)=\sum_{i=1}^{r}v^{(i)}\otimes u^{(i)}\otimes
f_{i}(x)$ and such that
\begin{eqnarray}
& &[\D\otimes 1,\S(x)]=-\frac{d}{dx}\S(x),\\
& &\S(z)(Y(x)\otimes 1) =(Y(x)\otimes 1)\S^{23}(z)\S^{13}(z+x).
\end{eqnarray}

\br{rnondegeneratcy} {\em A nonlocal vertex algebra $V$ is said to
be {\em nondegenerate} (see \cite{ek}) if for every positive integer
$n$, the linear map
$$Z_{n}: \C((x_{1}))\cdots ((x_{n}))\otimes V^{\otimes n}
\rightarrow V((x_{1}))\cdots ((x_{n}))$$ defined by
$$Z_{n}(f\otimes v^{(1)}\otimes \cdots \otimes v^{(n)})
=fY(v^{(1)},x_{1})\cdots Y(v^{(n)},x_{n}){\bf 1}$$ is injective. It
was proved in \cite{li-qva2} that if $V$ is of countable dimension
(over $\C$) and if $V$ as a (left) $V$-module is irreducible, then
$V$ is nondegenerate.} \er

\br{rwqva-qva} {\em A variation of a theorem of Etingof-Kazhdan
(\cite{ek}, Proposition 1.11), formulated in \cite{li-qva1}, is that
if a weak quantum vertex algebra $V$ is nondegenerate, then
$\S$-locality defines a unitary rational quantum Yang-Baxter
operator $\S(x)$ on $V$  such that $(V,\S(x))$ is a quantum vertex
algebra and this $\S(x)$ is the unique quantum Yang-Baxter operator
to make $V$ a quantum vertex algebra. In view of this, the term ``a
nondegenerate quantum vertex algebra,'' or ``an irreducible quantum
vertex algebra'' (without reference to a quantum Yang-Baxter
operator) is unambiguous. } \er

We shall need the following result of \cite{ll} (Proposition 4.5.7,
which can be extended for nonlocal vertex algebras with the same
proof):

\bl{lbasic-property-nva1} Let $V$ be a nonlocal vertex algebra, let
$W$ be a $V$-module, and let $u,v\in V,\; p,q\in \Z,\; w\in W$. Then
there exist nonnegative integers $l$ and $m$ such that
\begin{eqnarray}
u_{p}v_{q}w=\sum_{i=0}^{m}\sum_{j=0}^{l}\binom{p-l}{i}
\binom{l}{j}(u_{p-l-i+j}v)_{q+l+i-j}w.
\end{eqnarray}
\el

We shall also need the following analogue:

\bl{lbasic-property-nva} Let $V$ be a nonlocal vertex algebra and
let $u,v,w\in V,\; m,n\in \Z$. There exist nonnegative integers $l$
and $k$ such that
\begin{eqnarray}
(u_{m}v)_{n}w=\sum_{i=0}^{k}\sum_{j\ge 0}\binom{-l}{i}
\binom{m+i}{j}(-1)^{j}u_{m+l+i-j}v_{n-l-i+j}w.
\end{eqnarray}
\el

\begin{proof} {}From definition, there exists a nonnegative
integer $l$ such that
$$(x_{0}+x_{2})^{l}Y(Y(u,x_{0})v,x_{2})w
=(x_{0}+x_{2})^{l}Y(u,x_{0}+x_{2})Y(v,x_{2})w.$$
We have
\begin{eqnarray*}
(u_{m}v)_{n}w&=&\Res_{x_{0}}\Res_{x_{2}}x_{0}^{m}x_{2}^{n}
Y(Y(u,x_{0})v,x_{2})w\\
&=&\Res_{x_{0}}\Res_{x_{2}}x_{0}^{m}x_{2}^{n}
(x_{2}+x_{0})^{-l}\left((x_{2}+x_{0})^{l}Y(Y(u,x_{0})v,x_{2})w\right)\\
&=&\Res_{x_{0}}\Res_{x_{2}}\sum_{i\ge 0}\binom{-l}{i}x_{0}^{m+i}x_{2}^{n-l-i}
\left((x_{2}+x_{0})^{l}Y(Y(u,x_{0})v,x_{2})w\right).
\end{eqnarray*}
Let $k$ be a nonnegative integer such that
$x_{0}^{m+k}Y(u,x_{0})v\in V[[x_{0}]]$. Then
\begin{eqnarray*}
& &(u_{m}v)_{n}w\\
&=&\Res_{x_{0}}\Res_{x_{2}}\sum_{i=0}^{k}\binom{-l}{i}x_{0}^{m+i}x_{2}^{n-l-i}
\left((x_{2}+x_{0})^{l}Y(Y(u,x_{0})v,x_{2})w\right)\\
&=&\Res_{x_{0}}\Res_{x_{2}}\sum_{i=0}^{k}\binom{-l}{i}x_{0}^{m+i}x_{2}^{n-l-i}
\left((x_{0}+x_{2})^{l}Y(u,x_{0}+x_{2})Y(v,x_{2})w\right)\\
&=&\Res_{x_{1}}\Res_{x_{2}}\sum_{i=0}^{k}\binom{-l}{i}(x_{1}-x_{2})^{m+i}x_{2}^{n-l-i}
x_{1}^{l}Y(u,x_{1})Y(v,x_{2})w\\
&=&\Res_{x_{1}}\Res_{x_{2}} \sum_{i=0}^{k}\binom{-l}{i}\sum_{j\ge
0}\binom{m+i}{j}(-1)^{j}
x_{1}^{m+l+i-j}x_{2}^{n-l-i+j}Y(u,x_{1})Y(v,x_{2})w\\
&=&\sum_{i=0}^{k}\binom{-l}{i}
\sum_{j\ge 0}\binom{m+i}{j}(-1)^{j}u_{m+l+i-j}v_{n-l-i+j}w,
\end{eqnarray*}
as desired.
\end{proof}

Let $V$ be a nonlocal vertex algebra and let $T$ be a generating
subset of $V$ in the sense that $V$ is the smallest nonlocal vertex
subalgebra that contains $T$. Then (see \cite{bk}, \cite{li-g1}) $V$
is linearly spanned by the vectors
$$u^{(1)}_{m_{1}}\cdots u^{(r)}_{m_{r}}{\bf 1}$$
for $r\ge 0,\; u^{(1)},\dots,u^{(r)}\in T,\; m_{1},\dots,m_{r}\in
\Z$.

\br{rsimple-facts} {\rm  Here, we mention a simple fact which is
straightforward to prove. Let $V$ be a nonlocal vertex algebra with
a generating subset $T$ and let $K$ be any nonlocal vertex algebra.
For any map $f$ from $T$ to $K$, $f$ extends to at most one
nonlocal-vertex-algebra homomorphism from $V$ to $K$. On the other
hand, if $\psi$ is a linear map from $V$ to $K$ such that
$$\psi(Y(a,x)v)=Y(\psi(a),x)\psi(v)\ \ \ \mbox{ for }a\in T,\; v\in V,$$
then $\psi$ is a nonlocal-vertex-algebra homomorphism.} \er

In the following we shall associate two increasing filtrations of
$V$ to each generating subset.

\bl{lfiltration-E} Let $V$ be a nonlocal vertex algebra and let $T$
be a generating subset. For $n\in \N$, set
\begin{eqnarray}
E_{n}={\rm span} \{ u^{(1)}_{m_{1}}\cdots u^{(r)}_{m_{r}}{\bf 1}\;|\;
0\le r\le n,\; u^{(1)},\dots,u^{(r)}\in T,\; m_{1},\dots,m_{r}\in \Z\}.
\end{eqnarray}
Then $E_{n}\subset E_{n+1}$ for $n\in \N$, $E_{0}=\C {\bf 1}$, and
$\cup_{n\in \N}E_{n}=V$. Furthermore,
\begin{eqnarray}\label{e-filtration-t1}
a_{k}E_{n}\subset E_{m+n}\;\;\;\mbox{ for }a\in E_{m},\; k\in \Z,\;
m,n\in \N.
\end{eqnarray}
\el

\begin{proof} {}From definition, we have
$E_{n}\subset E_{n+1}$ for $n\in \N$ and $E_{0}=\C {\bf 1}$. As $T$
generates $V$, we get $V=\cup_{n\in \N}E_{n}$. With  $E_{0}=\C {\bf
1}$, we see that (\ref{e-filtration-t1}) holds for $m=0$. Assume
$m\ge 1$. {}From definition we have
$$E_{m}={\rm span}\{ u_{q}E_{m-1}\;|\; u\in T,\; q\in \Z\}.$$
It follows from induction on $m$ and Lemma \ref{lbasic-property-nva}
that (\ref{e-filtration-t1}) holds for all $m\in \N$.
\end{proof}

\bd{dvertex-graded-algebra} {\em Let $G$ be a group with identity
element $e$. A {\em nonlocal vertex $G$-graded algebra} is a
nonlocal vertex algebra $V$ equipped with a $G$-grading
$V=\coprod_{g\in G}V[g]$ such that ${\bf 1}\in V[e]$ and $Y(u,x)v\in
V[gh]((x))$ for $u\in V[g],\; v\in V[h]$ with $g,h\in G$. } \ed

 \br{rfiltration-E} {\em Let
$V$ be a nonlocal vertex algebra and let $E=\{ E_{n}\}_{n\in \Z}$ be
an increasing subspace-filtration of $V$ with ${\bf 1}\in E_{0}$,
satisfying the condition that
\begin{eqnarray*}
a_{k}E_{n}\subset E_{m+n}\;\;\;\mbox{ for }a\in E_{m},\; m,n,k\in \Z
\end{eqnarray*}
(cf. (\ref{e-filtration-t1})). Form the graded vector space
\begin{eqnarray}
{\rm Gr}_{E}(V)=\coprod_{n\in \Z}(E_{n}/E_{n-1}).
\end{eqnarray}
It is straightforward to show (cf. \cite{li-qva2}, Lemma 3.13) that
${\rm Gr}_{E}(V)$ is a nonlocal vertex algebra with ${\bf
1}+E_{-1}\in E_{0}/E_{-1}$ as the vacuum vector and with
\begin{eqnarray}
(a+E_{m-1})_{k}(b+E_{n-1})=a_{k}b+E_{m+n-1}\in E_{m+n}/E_{m+n-1}
\end{eqnarray}
for $a\in E_{m},\; b\in E_{n}$, $m,n,k\in \Z$. In fact, ${\rm
Gr}_{E}(V)$ is a nonlocal vertex $\Z$-graded algebra.   Furthermore,
if $E$ is the filtration associated with a generating subset $T$ of
$V$, then $\{ u+E_{0}\;|\; u\in T\}$ $(\subset E_{1}/E_{0})$ is a
generating subset of ${\rm Gr}_{E}(V)$. A result that was obtained
in \cite{li-qva2} is that if ${\rm Gr}_{E}(V)$ is nondegenerate,
then $V$ is nondegenerate.} \er

With a generating subset $T$ of $V$ we have another increasing
filtration.

\bl{linduction1} Let $V$ be a nonlocal vertex algebra and let $T$ be
a generating subset. For $n\in \Z$, let $F_{n}$ be the linear span
of the vectors
$$u^{(1)}_{m_{1}}\cdots u^{(r)}_{m_{r}}{\bf 1}$$
for $r\ge 1$ if $n<0$, for $r\ge 0$ if $n\ge 0$, and for
$u^{(1)},\dots,u^{(r)}\in T,\; m_{1},\dots,m_{r}\in \Z$ with
$m_{1}+\cdots +m_{r}\ge -n$. Then $\{F_{n}\}_{n\in \Z}$ is an
increasing filtration of $V$ such that ${\bf 1}\in F_{0}$ and
\begin{eqnarray}\label{e-filtration-t2}
a_{m}F_{n}\subset F_{k+n-m-1} \;\;\;\mbox{ for }a\in F_{k},\;
k,m,n\in\Z.
\end{eqnarray}
\el

\begin{proof} It is clear that ${\bf 1}\in F_{0}$ and
$\{F_{n}\}_{n\in \Z}$ is an increasing filtration of $V$. For $u\in
T,\; m,n\in \Z$, from definition we have
$$u_{m}F_{n}\subset F_{n-m}.$$
It follows from induction on $r$
and Lemma \ref{lbasic-property-nva} that
for any nonnegative integer $r$,
for any $a$ of the form $u^{(1)}_{m_{1}}\cdots u^{(r)}_{m_{r}}{\bf 1}$
with $u^{(1)},\dots,u^{(r)}\in T,\; m_{1},\dots,m_{r}\in \Z$,
and for any $m\in \Z$ we have
$$a_{m}F_{n}\subset F_{n-m-m_{1}-\cdots -m_{r}-1}.$$
Then (\ref{e-filtration-t2}) is clear.
\end{proof}

A {\em $\Z$-graded nonlocal vertex algebra} is a nonlocal vertex
algebra $U$ equipped with a $\Z$-grading $U=\coprod_{n\in
\Z}U_{(n)}$ such that ${\bf 1}\in U_{(0)}$ and
$$u_{k}v\in U_{(m+n-k-1)}$$
for all $u\in U_{(m)},\; v\in U_{(n)},\; m,n,k\in \Z$.

 \bp{pZ-gradedva} Let $F=\{
F_{n}\}_{n\in \Z}$ be an increasing filtration of $V$ with ${\bf
1}\in F_{0}$ satisfying (\ref{e-filtration-t2}). Then the associated
graded vector space ${\rm Gr}_{F}(V)=\coprod_{n\in
\Z}(F_{n}/F_{n-1})$ is a $\Z$-graded nonlocal vertex algebra with
${\bf 1}+F_{-1}$ as the vacuum vector, where
$$(a+F_{m-1})_{k}(b+F_{n-1})=a_{k}b+F_{m+n-k-2}$$
for $a\in F_{m},\; b\in F_{n},\; m,n,k\in \Z$.  Furthermore, if $F$
is the filtration associated with a generating subset $T$ of $V$,
then $\{u+F_{0}\;|\; u\in T\}\subset F_{1}/F_{0}$ is a generating
subset of
 ${\rm Gr}_{F}(V)$. \ep

\begin{proof} It is easy to see that the axioms that involve the vacuum
vector, namely the vacuum and creation properties, hold. For weak
associativity, we need to show that for $a\in F_{m},\; b\in F_{n},\;
c\in F_{k}$ with $m,n,k\in \Z$, there exists a nonnegative integer
$l$ such that
\begin{eqnarray}\label{ewa-barabc}
(x_{0}+x_{2})^{l}Y(\bar{a},x_{0}+x_{2})Y(\bar{b},x_{2})\bar{c}
=(x_{0}+x_{2})^{l}Y(Y(\bar{a},x_{0})\bar{b},x_{2})\bar{c},
\end{eqnarray}
where $\bar{a}=a+F_{m-1},\; \bar{b}=b+F_{n-1},\; \bar{c}=c+F_{k-1}$.
Let $l$ be a nonnegative integer such that
$$(x_{0}+x_{2})^{l}Y(a,x_{0}+x_{2})Y(b,x_{2})c
=(x_{0}+x_{2})^{l}Y(Y(a,x_{0})b,x_{2})c,$$ which is expanded as
$$\sum_{p,q\in \Z}\sum_{i\ge 0}\binom{l-p-1}{i}a_{p}b_{q}c
x_{0}^{l-p-1-i}x_{2}^{i-q-1} =\sum_{p',q'\in \Z}\sum_{i'\ge
0}\binom{l}{i'}(a_{p'}b)_{q'}cx_{0}^{l-i'-p'-1}x_{2}^{i'-q'-1}.$$
Let $r,s\in \Z$. Extracting the coefficients of
$x_{0}^{r}x_{2}^{s}$ from both sides we get
\begin{eqnarray}\label{ecomponents-rs}
\sum_{i\ge 0}\binom{r+i}{i}a_{l-r-1-i}b_{i-1-s}c
 =\sum_{i'\ge
0}\binom{l}{i'}(a_{l-i'-r-1}b)_{i'-s-1}c.
\end{eqnarray}
With $a\in F_{m},\; b\in F_{n},\; c\in F_{k}$, by definition we have
$$\bar{b}_{i-1-s}\bar{c}=b_{i-1-s}c+F_{n+k+s-i-1},$$
and furthermore,
\begin{eqnarray}\label{ebarabc-formula}
\bar{a}_{l-r-1-i}\bar{b}_{i-1-s}\bar{c}
=a_{l-r-1-i}b_{i-1-s}c+F_{m+n+k+r+s-l-1}.
\end{eqnarray}
On the other hand, we have
$$ \bar{a}_{l-i'-r-1}\bar{b}=a_{l-i'-r-1}b+F_{m+n+r+i'-l-1}$$ and
\begin{eqnarray}\label{ebarabc-iterate}
(\bar{a}_{l-i'-r-1}\bar{b})_{i'-s-1}\bar{c}
=(a_{l-i'-r-1}b)_{i'-s-1}c+F_{m+n+k+r+s-l-1}.
\end{eqnarray}
Then by (\ref{ecomponents-rs}) we have
$$\sum_{i\ge 0}\binom{r+i}{i}\bar{a}_{l-r-1-i}\bar{b}_{i-1-s}\bar{c}
 =\sum_{i'\ge
0}\binom{l}{i'}(\bar{a}_{l-i'-r-1}\bar{b})_{i'-s-1}\bar{c}.$$
Multiplying both sides by $x_{0}^{r}x_{2}^{s}$, summing up over all
$r,s\in \Z$, and then changing indices we get
$$\sum_{p,q\in \Z}\sum_{i\ge
0}\binom{l-p-1}{i}\bar{a}_{p}\bar{b}_{q}\bar{c}
x_{0}^{l-p-1-i}x_{2}^{i-q-1} =\sum_{p',q'\in \Z}\sum_{i'\ge
0}\binom{l}{i'}(\bar{a}_{p'}\bar{b})_{q'}\bar{c}x_{0}^{l-i'-p'-1}x_{2}^{i'-q'-1},$$
which gives (\ref{ewa-barabc}).
 This proves that ${\rm Gr}_{F}(V)$ is a nonlocal vertex
algebra. It follows from definition that ${\rm Gr}_{F}(V)$ is a
$\Z$-graded nonlocal vertex algebra.

For the last assertion, from definition we have $T\subset F_{1}$ (as
$u=u_{-1}{\bf 1}$ for $u\in T$). For any $n\in \Z$, consider a
typical vector
$$a=u^{(1)}_{m_{1}}\cdots u^{(r)}_{m_{r}}{\bf 1}\in
F_{n},$$ where $r\ge 1,\; u^{(1)},\dots,u^{(r)}\in T,\;
m_{1},\dots,m_{r}\in \Z$ with $m_{1}+\cdots +m_{r}\ge -n$. If
$m_{1}+\cdots +m_{r}>-n$, we have $a\in F_{n-1}$, so that
$a+F_{n-1}=0$ in ${\rm Gr}_{F}(V)$. Assume that $m_{1}+\cdots
+m_{r}=-n$. Set $\bar{\bf 1}={\bf 1}+F_{-1}$. Note that
$u_{m}F_{k}\subset F_{k-m}$ for $u\in T,\; m,k\in \Z$. From
definition we have
$$a+F_{n-1}=u^{(1)}_{m_{1}}\cdots u^{(r)}_{m_{r}}{\bf
1}+F_{-m_{1}-\cdots -m_{r}-1} =\bar{u}^{(1)}_{m_{1}}\cdots
\bar{u}^{(r)}_{m_{r}}\bar{\bf 1}.$$ Then it follows that ${\rm
Gr}_{F}(V)$ is generated by $\{ u+F_{0}\;|\; u\in T\}$.
\end{proof}

\bp{p-filtration-irred} Let $V$ be a nonlocal vertex algebra. Assume
that there exists a lower-truncated increasing filtration $F=\{
F_{n}\}_{n\in \Z}$ of $V$ with ${\bf 1}\in F_{0}$, satisfying either
(\ref{e-filtration-t1}) or (\ref{e-filtration-t2}), such that the
left adjoint module for ${\rm Gr}_{F}(V)$ is graded-irreducible.
Then the left adjoint module for $V$ is irreducible. \ep

\begin{proof} Let $U$ be a nonzero $V$-submodule of $V$.
We have an increasing filtration $F_{U}=\{ F_{n}\cap U\}_{n\in \Z}$
of $U$. For each $n\in \Z$, the natural map {}from $F_{n}\cap U$ to
$F_{n}/F_{n-1}$ gives rise to an embedding of $(F_{n}\cap
U)/(F_{n-1}\cap U)$ into $F_{n}/F_{n-1}$. Then we can view ${\rm
Gr}_{F_{U}}(U)$ as a subspace of ${\rm Gr}_{F}(V)$. As $U$ is a
$V$-submodule, ${\rm Gr}_{F_{U}}(U)$ is a graded ${\rm
Gr}_{F}(V)$-submodule. Since ${\rm Gr}_{F}(V)$ is
graded-irreducible, we have either ${\rm Gr}_{F_{U}}(U)=0$ or ${\rm
Gr}_{F_{U}}(U)={\rm Gr}_{F}(V)$.

Assume  ${\rm Gr}_{F_{U}}(U)=0$. Then $F_{n}\cap U=F_{n-1}\cap U$
for all $n\in \Z$. For any $k\in \Z$, we have $F_{k}\cap
U=\cup_{n\ge k}(F_{n}\cap U)=U$, which implies $U\subset F_{k}$.
Thus $U\subset \cap_{k\in \Z}F_{k}$. As $F$ is lower-truncated, we
have $\cap_{n\in \Z} F_{n}=0$. This contradicts that $U\ne 0$.
Therefore, we must have ${\rm Gr}_{F_{U}}(U)={\rm Gr}_{F}(V)$. Then
$(F_{n+1}\cap U)/(F_{n}\cap U)=F_{n+1}/F_{n}$ for $n\in \Z$. That
is, $F_{n+1}=(F_{n+1}\cap U)+F_{n}\subset U+F_{n}$ for all $n\in
\Z$. {}From this we get
$$F_{n+i}\subset U+F_{n}\ \ \ \mbox{ for }i\ge 0.$$
For every $n\in \Z$, as $\cup_{i\ge 0}F_{n+i}=V$, we have
$V=U+F_{n}$. Then $V=U$ because $F_{n}=0$ for some $n\in \Z$. This
proves that $V$ is irreducible.
\end{proof}

\br{rirred} {\em Let $V=\coprod_{n\in \Z}V_{(n)}$ be a lower
truncated $\Z$-graded nonlocal vertex algebra. For $n\in \Z$, set
$F_{n}=\coprod_{m\le n}V_{(m)}$. Then $\{ F_{n}\}$ is an increasing
filtration satisfying (\ref{e-filtration-t2}). Furthermore,
$Gr_{F}(V)\simeq V$ as a nonlocal vertex algebra. It follows from
Proposition \ref{p-filtration-irred} that if $V$ is graded
irreducible, then $V$ is irreducible.} \er

The following is straightforward:

\bl{ltwisted} Let $V$ be a nonlocal vertex $G$-graded algebra with
$G$ an abelian group and let $\varepsilon: G\times G\rightarrow
\C^{\times}$ be a normalized $\C^{\times}$-valued $2$-cocycle of $G$
in the sense that
\begin{eqnarray*}
& &\varepsilon(\alpha,0)=\varepsilon(0,\alpha)=1,\\
&& \varepsilon(\alpha,\beta+\gamma)\varepsilon(\beta,\gamma)
=\varepsilon(\alpha,\beta)\varepsilon(\alpha+\beta,\gamma) \ \ \
\mbox{ for }\alpha,\beta,\gamma\in G.
\end{eqnarray*}
 Define a linear map $Y_{\varepsilon}: V\rightarrow (\End V)[[x,x^{-1}]]$ by
 $$Y_{\varepsilon}(u,x)v=\varepsilon(\alpha,\beta)Y(u,x)v\;\;\;\mbox{
for }u\in V[\alpha],\; v\in V[\beta],\;\alpha,\beta\in G.$$
 Then
$(V,Y_{\varepsilon}, {\bf 1})$ carries the structure of a nonlocal
vertex $G$-graded algebra. Furthermore, $(V,Y_{\varepsilon}, {\bf
1})$ is nondegenerate if and only if $V$ is nondegenerate, and
$(V,Y_{\varepsilon}, {\bf 1})$ is $G$-graded irreducible if and only
if $V$ is $G$-graded irreducible. \el

The following is a vertex analogue of a well known construction of
associative algebras:

\bp{psmash} Let $G$ be a group, let $U$ be a nonlocal vertex algebra
on which $G$ acts by automorphisms, and let $V=\coprod_{g\in G}V[g]$
be a nonlocal vertex $G$-graded algebra. Define a linear map
$$Y_{\sharp}: (U\otimes V)\rightarrow (\End (U\otimes V))[[x,x^{-1}]]$$
by
\begin{eqnarray}
Y_{\sharp}(u\otimes v,x)(u'\otimes v')=Y(u,x)g(u')\otimes Y(v,x)v'
\end{eqnarray}
for $u,u'\in U,\; v\in V[g],\; g\in G,\; v'\in V$. Then $(U\otimes
V, Y_{\sharp},{\bf 1}\otimes {\bf 1})$ carries the structure of a
nonlocal vertex algebra which we denote by $U\sharp_{G}V$.
Furthermore, $U$ and $V$ are subalgebras of $U\sharp_{G}V$ and
\begin{eqnarray}
Y_{\sharp}(v,x_{1})Y_{\sharp}(u,x_{2})
=Y_{\sharp}(g(u),x_{2})Y_{\sharp}(v,x_{1})
\end{eqnarray}
for $u\in U,\; v\in V[g]$ with $g\in G$. Moreover, $U\sharp_{G}V$ is
a weak quantum vertex algebra if both $U$ and $V$ are weak quantum
vertex algebras. \ep

\begin{proof} First, for $u,u'\in U,\; v\in V[g],\;v'\in V$ with $g\in G$ we have
$$Y_{\sharp}(u\otimes v,x)(u'\otimes v')=Y(u,x)g(u')\otimes Y(v,x)v'
\in (U\otimes V)((x)).$$  Second, $Y_{\sharp}({\bf 1}\otimes {\bf
1},x)(u\otimes v)=u\otimes v$ and $Y_{\sharp}(u\otimes v,x)({\bf
1}\otimes {\bf 1})=Y(u,x){\bf 1}\otimes Y(v,x){\bf 1}.$ Let
$u,u',u''\in U,\; v\in V[g],\;v'\in V[h],\;v''\in V$ with $g,h\in
G$. Then
\begin{eqnarray*}
& &Y_{\sharp}(u\otimes v,x_{0}+x_{2})Y_{\sharp}(u'\otimes
v',x_{2})(u''\otimes v'')\\
&=&Y(u,x_{0}+x_{2})gY(u',x_{2})hu''\otimes
Y(v,x_{0}+x_{2})Y(v',x_{2})v''\\
&=&Y(u,x_{0}+x_{2})Y(gu',x_{2})ghu''\otimes
Y(v,x_{0}+x_{2})Y(v',x_{2})v'',
\end{eqnarray*}
and
\begin{eqnarray*}
& &Y_{\sharp}(Y_{\sharp}(u\otimes v,x_{0})(u'\otimes
v'),x_{2})(u''\otimes v'')\\
&=&Y_{\sharp}(Y(u,x_{0})gu'\otimes Y(v,x_{0})v',x_{2})(u''\otimes
v'')\\
&=&Y(Y(u,x_{0})gu',x_{2})(gh)u''\otimes Y(Y(v,x_{0})v',x_{2})v'',
\end{eqnarray*}
as $Y(v,x_{0})v'\in V[gh]((x_{0}))$. Then weak associativity follows
immediately. This proves that $(U\otimes V,Y_{\sharp},{\bf 1}\otimes
{\bf 1})$ is a nonlocal vertex algebra. Recall that ${\bf 1}\in
V[e]$ and that $e$ acts on $U$ as identity. It is easy to see that
$U$, identified with $U\otimes {\bf 1}$, and $V$, identified with
${\bf 1}\otimes V$, are subalgebras of $U\sharp_{G} V$.

For the last assertion, let $u,u'\in U,\; v\in V[g],\; v'\in V$.
Then
\begin{eqnarray*}
&&Y_{\sharp}({\bf 1}\otimes v,x_{1})Y_{\sharp}(u\otimes {\bf
1},x_{2})(u'\otimes v')\\
&=&gY(u,x_{2})eu'\otimes Y(v,x_{1})v'\\
&=&Y(gu,x_{2})gu'\otimes Y(v,x_{1})v'\\
&=& Y_{\sharp}(gu\otimes {\bf 1},x_{2})Y_{\sharp}({\bf 1}\otimes
v,x_{1})(u'\otimes v').
\end{eqnarray*}

For $u\in U,\; v\in V$, we have
$$Y_{\sharp}(u\otimes {\bf 1},x)({\bf 1}\otimes v)=Y(u,x){\bf
1}\otimes v.$$ It follows that $U\sharp_{G}V$ is generated by $U$
and $V$. Then it is clear that  $U\sharp_{G}V$ is a weak quantum
vertex algebra if both $U$ and $V$ are weak quantum vertex algebras.
\end{proof}

The smash product nonlocal vertex algebra $U\sharp_{G}V$ has the
following universal property:

\bp{psmash-universal} Let $G,U,V$ be given as in Proposition
\ref{psmash}. Suppose that $K$ is a nonlocal vertex algebra and
$\psi: U\rightarrow K$ and $\phi: V\rightarrow K$ are homomorphisms
of nonlocal vertex algebras such that
\begin{eqnarray}\label{e2.24}
Y(\phi(v),x_{1})Y(\psi(u),x_{2})=Y(\psi(gu),x_{2})Y(\phi(v),x_{1})
\end{eqnarray}
for $u\in U,\; v\in V[g]$ with $g\in G$. Then there exists a unique
nonlocal vertex algebra homomorphism $f: U\sharp_{G}V\rightarrow K$,
extending both $\psi$ and $\phi$. \ep

\begin{proof} As $U\sharp_{G}V$ is generated by $U$ and
$V$, the uniqueness follows immediately. For the existence, we
define a linear map $\theta: U\sharp_{G}V\rightarrow K$ by
$$\theta(u\otimes v)=\psi(u)_{-1}\phi(v)\ \ \ \mbox{ for }u\in U,\; v\in V.$$
It is clear that $\theta$ extends both $\psi$ and $\phi$. It remains
to prove that $\theta$ is a homomorphism of nonlocal vertex
algebras. For $u\in U,\; v\in V[g]$, by (\ref{e2.24}) we have
$$Y(\psi(u),x_{2})Y(\phi(v),x_{1}){\bf
1}=Y(\phi(v),x_{1})Y(\psi(g^{-1}u),x_{2}){\bf 1},$$ which implies
that $Y(\psi(u),x)\phi(v)\in K[[x]]$. For any $u\in U,\; v,v'\in V$,
there exists a nonnegative integer $l$ such that
$$ (x_{0}+x_{2})^{l}Y(\psi(u),x_{0}+x_{2})Y(\phi(v),x_{2})\phi(v')
= (x_{0}+x_{2})^{l}Y(Y(\psi(u),x_{0})\phi(v),x_{2})\phi(v'). $$ As
$Y(\psi(u),x)Y(\phi(v),x_{2})\phi(v')$ involves only nonnegative
powers of $x$, in fact we have
$$
Y(\psi(u),x_{0}+x_{2})Y(\phi(v),x_{2})\phi(v') =
Y(Y(\psi(u),x_{0})\phi(v),x_{2})\phi(v'). $$
 Furthermore, we have
 \begin{eqnarray*}
 \theta(Y_{\sharp}(u\otimes
{\bf 1},x)(u'\otimes v'))&=&\theta(Y(u,x)u'\otimes v')\\
&=&\lim_{x_{2}\rightarrow
0}Y(\psi(Y(u,x)u'),x_{2})\phi(v')\\
&=&\lim_{x_{2}\rightarrow 0}Y(Y(\psi(u),x)\psi(u'),x_{2})\phi(v')\\
&=&\lim_{x_{2}\rightarrow 0}Y(\psi(u),x+x_{2})Y(\psi(u'),x_{2})\phi(v')\\
&=&Y(\psi(u),x)\theta(u'\otimes v')
\end{eqnarray*}
and
\begin{eqnarray*}
 \theta(Y_{\sharp}({\bf 1}\otimes
v,x)(u'\otimes v'))&=&\theta(gu'\otimes Y(v,x)v')\\
&=&\lim_{x_{2}\rightarrow
0}Y(\psi(gu'),x_{2})\phi(Y(v,x)v')\\
&=&\lim_{x_{2}\rightarrow 0}Y(\psi(gu'),x_{2})Y(\phi(v),x)\phi(v')\\
&=&\lim_{x_{2}\rightarrow 0}Y(\phi(v),x+x_{2})Y(\psi(u'),x_{2})\phi(v')\\
&=&Y(\phi(v),x)\theta(u'\otimes v').
\end{eqnarray*}
As $U$ and $V$ generate $U\sharp_{G}V$, it follows that $\theta$ is
a homomorphism.
\end{proof}

We continue to study $U\sharp_{G}V$-modules.

\bp{pproduct-module} Let $E$ be a $U$-module on which $G$ acts such
that
$$g(Y(u,x)w)=Y(gu,x)gw\;\;\;\mbox{ for }g\in G,\; u\in U,\; w\in E,$$
and let $F$ be any $V$-module. Then $E\otimes F$ is a $U\sharp_{G}
V$-module with
\begin{eqnarray}
Y(u\otimes v,x)(w\otimes f)=Y(u,x)gw\otimes Y(v,x)f
\end{eqnarray}
for $u\in U,\; v\in V[g],\;g\in G,\; w\in E,\; f\in F$. We denote
this module by $E\sharp F$. Furthermore, if the $U$-module $E$ and
the $V$-module $F$ are irreducible and if $E$ is of countable
dimension (over $\C$), then $E\sharp F$ is an irreducible
$U\sharp_{G} V$-module. \ep

\begin{proof} From the first part of the proof of Proposition \ref{psmash},
we see that this indeed defines a $U\sharp_{G}V$-module structure on
$E\otimes F$. For proving irreducibility, let $W$ be any nonzero
submodule of $E\otimes F$.  Notice that $U$ acts on the first factor
of $E\otimes F$. Since $E$ is an irreducible $U$-module of countable
dimension (over $\C$), the Schur lemma holds and $\End_{U}E=\C$.
Using Jacobson's density theorem on $E$, we get $W=E\otimes K$ for
some subspace $K$ of $F$. {}From the definition of the action of $V$
on $E\otimes F$, we see that $K$ is a $V$-submodule of $F$. With $F$
an irreducible $V$-module we have $K=F$, proving $W=E\otimes F$.
Thus $E\otimes F$ is an irreducible $U\sharp_{G} V$-module.
\end{proof}

\bc{cnondegeneracy} Assume that both $U$ and $V$ are irreducible
nonlocal vertex algebras and that $U$ is of countable dimension
(over $\C$). Then $U\sharp_{G} V$ is an irreducible nonlocal vertex
algebra. \ec

\bp{pmodule-class} Let $W$ be a $U\sharp_{G}V$-module. Assume that
$W$ contains an irreducible $U$-submodule $E$ of countable dimension
(over $\C$) on which $G$ acts such that
$$g(Y(u,x)w)=Y(gu,x)gw\;\;\;\mbox{ for }g\in G,\; u\in U,\; w\in E$$
and assume that $W$ is generated by $E$. Then $W$ is isomorphic to
$E\otimes F$ for some $V$-module $F$. Furthermore, if $W$ is
irreducible, then $F$ is an irreducible $V$-module.
 \ep

\begin{proof} First we prove that $W$ as a $U$-module is a
direct sum of irreducible modules isomorphic to $E$.
 Let $v\in V[g],\; n\in \Z$ with $g\in G$. It
is clear that $v_{n}E$ is a $U$-submodule of $W$ and the linear map
from $E$ to $v_{n}E$, sending $w$ to $v_{n}g^{-1}w$, is a
$U$-homomorphism. Thus $v_{n}E$, if not zero, is an irreducible
$U$-submodule of $W$, isomorphic to $E$. Let $W'$ be the sum of the
$U$-submodules $v_{n}E$ for $v\in V[g],\; n\in \Z$ with $g\in G$. It
follows from Lemma \ref{lbasic-property-nva1} that $W'$ is a
$V$-submodule of $W$. As $U$ and $V$ generate $U\sharp_{G}V$, it
follows that $W'$ is a $U\sharp_{G}V$-submodule. Thus $W'=W$ as $W$
is generated by $E$. This proves that $W$ as a $U$-module is a
direct sum of copies of $E$.

With $E$ an irreducible $U$-module of countable dimension (over
$\C$), we have $W=E\otimes \Hom _{U}(E,W)$. For $v\in V[g],\; g\in
G,\; f\in \Hom_{U}(E,W)$, we define $Y(v,x)f\in
(\Hom_{\C}(E,W))[[x,x^{-1}]]$ by
$$(Y(v,x)f)(w)=Y(v,x)f(g^{-1}w)\;\;\;\mbox{ for }w\in E.$$
For $u\in U,\; w\in W$, we have
\begin{eqnarray*}
(Y(v,x)f)(Y(u,x_{0})w)&=&Y(v,x)f(g^{-1}Y(u,x_{0})w)\\
&=&Y(v,x)f(Y(g^{-1}u,x_{0})g^{-1}w)\\
&=&Y(v,x)Y(g^{-1}u,x_{0})f(g^{-1}w)\\
&=&Y(u,x_{0})Y(v,x)f(g^{-1}w)\\
&=&Y(u,x_{0})(Y(v,x)f)(w).
\end{eqnarray*}
This proves that $Y(v,x)f\in (\Hom_{U}(E,W))[[x,x^{-1}]]$. Let $0\ne
w\in E$ be arbitrarily fixed. Since $E$ is an irreducible
$U$-module, $g^{-1}w$ generates $E$. Let $l$ be a nonnegative
integer such that $x^{l}Y(v,x)f(g^{-1}w)\in W[[x]]$. It follows that
$$x^{l}(Y(v,x)f)\in (\Hom_{U}(E,W))[[x]].$$
That is, $Y(v,x)f\in (\Hom_{U}(E,W))((x))$. Furthermore, for $v\in
V[g],\; v'\in V[h]$ with $g,h\in G$ and for $w\in W$, we have
\begin{eqnarray*}
(Y(v,x_{1})Y(v',x_{2})f)(w)=Y(v,x_{1})(Y(v',x_{2})f)(g^{-1}w)
=Y(v,x_{1})Y(v',x_{2})f(h^{-1}g^{-1}w)
\end{eqnarray*}
and
$$(Y(Y(v,x_{0})v',x_{2})f)(w)=Y(Y(v,x_{0})v',x_{2})f((gh)^{-1}w).$$
Now weak associativity follows. This proves that we have a
$V$-module structure on $\Hom_{U}(E,W)$, so that $E\otimes \Hom
_{U}(E,W)$ is a $U\sharp_{G}V$-module by Proposition
\ref{pproduct-module} and we have $W=E\otimes \Hom _{U}(E,W)$ as a
$U\sharp_{G}V$-module. If $W$ is irreducible, $\Hom_{U}(E,W)$ must
be an irreducible $V$-module.
\end{proof}

\section{Associative algebras $\A_{\bf Q}$ and quantum vertex algebras
$V_{\bf Q}$}

In this section we associate an associative algebra $\A_{\bf Q}$ to
each ``skew'' complex matrix ${\bf Q}$ and we study its modules,
including vacuum modules. We show that on the universal vacuum
module $V_{\bf Q}$ there exists a canonical irreducible quantum
vertex algebra structure.

\bd{dquantum-torus} {\em Let $l$ be a positive integer and let ${\bf
Q}=(q_{ij})_{i,j=1}^{l}$ be a complex matrix such that
\begin{eqnarray}\label{eq-skewsymmetry}
q_{ij}q_{ji}=1\ \ \ \mbox{ for }1\le i,j\le l.
\end{eqnarray}
Define $\A_{\bf Q}$ to be the associative algebra with identity
(over $\C$) with generators
$$X_{i,n},\; Y_{i,n}\;\;\ \ (i=1,\dots,l,\; n\in \Z),$$
subject to relations
\begin{eqnarray}
& &X_{i,m}X_{j,n}=q_{ij}X_{j,n}X_{i,m},
\ \ \ \ Y_{i,m}Y_{j,n}=q_{ij}Y_{j,n}Y_{i,m},\nonumber\\
&&X_{i,m}Y_{j,n}-q_{ji}Y_{j,n}X_{i,m}=\delta_{i,j}\delta_{m+n+1,0}
\end{eqnarray}
for $i,j=1,\dots, l,\; m,n\in \Z$.} \ed

Let $\{e_{1},\dots,e_{l}\}$ denote the standard $\Z$-basis of
$\Z^{l}$. It is straightforward to see that $\A_{\bf Q}$ is a
$\Z^{l}$-graded algebra with the grading defined by
\begin{eqnarray}
\deg X_{i,m}=e_{i},\ \ \ \ \ \deg Y_{i,m}=-e_{i} \ \ \ \mbox{ for
}1\le i\le l,\; m\in \Z.
\end{eqnarray}
 Set
\begin{eqnarray*} &
&\A_{\bf Q}^{+}=\<X_{i,m},\; Y_{j,n}\;|\; i,j=1,\dots,l,\; m,n\ge
0\>,\\
& &\A_{\bf Q}^{-}=\<X_{i,m},\; Y_{j,n}\;|\;i,j=1,\dots,l,\; m,n<0\>,
\end{eqnarray*}
which are $\Z^{l}$-graded subalgebras of $\A_{\bf Q}$.

 Notice that if $q_{ij}=1$ for all
$i,j=1,\dots,l$, the algebra $\A_{\bf Q}$, which is isomorphic to
the universal enveloping algebra of an infinite-dimensional
Heisenberg Lie algebra,
 is a Weyl algebra, and that if $q_{ij}=-1$ for
all $i,j=1,\dots,l$, $\A_{\bf Q}$ is a Clifford algebra. In general,
we have $q_{ii}=\pm 1$ for $1\le i\le l$ as $q_{ii}q_{ii}=1$. Then,
for each $1\le i\le l$, the algebra $\A_{(q_{ii})}$ (associated with
the $1\times 1$ matrix $q_{ii}$) is either a Weyl algebra or a
Clifford algebra. In the following we shall prove that $\A_{\bf Q}$
is isomorphic to a certain cocycle-twist of the tensor product
algebra $\A_{(q_{11})}\otimes \cdots \otimes \A_{(q_{ll})}$.

For $1\le i\le l$, we temporarily denote the generators of
$\A_{(q_{ii})}$ by $\bar{X}_{i,m}$ and $\bar{Y}_{i,m}$ for $m\in
\Z$. Define a $\Z$-grading on $\A_{(q_{ii})}$ by
\begin{eqnarray}
\deg \bar{X}_{i,m}=1,\ \ \ \  \deg \bar{Y}_{i,m}=-1\ \ \ \mbox{ for
}m\in \Z,
\end{eqnarray}
making $\A_{(q_{ii})}$ a $\Z$-graded algebra, where we denote by
$\A_{(q_{ii})}[n]$ the degree-$n$ subspace for $n\in \Z$. Set
$$\A=\A_{(q_{11})}\otimes \A_{(q_{22})}\otimes \cdots \otimes
\A_{(q_{ll})},$$ the tensor product algebra, which is naturally a
$\Z^{l}$-graded algebra.  We define a group homomorphism
$\varepsilon: \Z^{l}\times \Z^{l}\rightarrow \C^{\times}$ by
\begin{eqnarray}\label{def-varepsilon}
\varepsilon (e_{i},e_{j})=\begin{cases}q_{ij}\ \ \ \mbox{ if }i> j\\
1\ \ \ \ \ \mbox{ if }i\le j\end{cases}
\end{eqnarray}
for $1\le i,j\le l$, recalling that $\{e_{1},\dots,e_{l}\}$ is the
standard $\Z$-basis of $\Z^{l}$. In particular, $\varepsilon$ is a
normalized $2$-cocycle of $\Z^{l}$. Denote by $\A^{\varepsilon}$ the
$\varepsilon$-twist of $\A$. That is, $\A^{\varepsilon}=\A$ as a
vector space and for $a\in \A[\alpha],\; b\in \A[\beta]$ with
$\alpha,\beta\in \Z^{l}$,
\begin{eqnarray}
a^{\varepsilon}b^{\varepsilon}=\varepsilon(\alpha,\beta)(ab)^{\varepsilon},
\end{eqnarray}
where for any $u\in \A$, we use $u^{\varepsilon}$ for $u$ viewed as
an element of $\A^{\varepsilon}$. Notice that for each $1\le i\le
l$, as $\varepsilon(e_{i},e_{i})=1$,
$\A_{(q_{ii})}^{\varepsilon}=\A_{(q_{ii})}$ as an algebra.

 \bp{ptwisted-algebra} The map $$\phi:
\{ \bar{X}_{i,m}^{\epsilon}, \bar{Y}_{i,m}^{\epsilon}\;|\; 1\le i\le
l,\; m\in \Z\}\subset \A^{\varepsilon} \rightarrow \A_{\bf Q};\ \
\bar{X}_{i,m}^{\epsilon}\mapsto X_{i,m},\ \
\bar{Y}_{i,m}^{\epsilon}\mapsto Y_{i,m}$$ extends to an algebra
isomorphism from  the $\varepsilon$-twist $\A^{\varepsilon}$ of the
$\Z^{l}$-graded algebra $\A$ to $\A_{\bf Q}$. Furthermore, for each
$1\le i\le l$, the assignment $$\bar{X}_{i,m}\mapsto  X_{i,m},\ \ \
\ \bar{Y}_{i,m}\mapsto Y_{i,m} \ \ \ \mbox{ for }m\in \Z$$ gives
rise to an algebra embedding of $\A_{(q_{ii})}$ into $\A_{\bf Q}$,
and the linear map
\begin{eqnarray}
\pi: \A_{(q_{11})}\otimes \cdots \otimes \A_{(q_{ll})}\rightarrow
\A_{\bf Q};\;\; (a_{1},\dots,a_{l})\mapsto a_{1}\cdots a_{l}
\end{eqnarray}
is a linear isomorphism preserving the $\Z^{l}$-gradings. \ep

\begin{proof} {}From the definition of $\varepsilon$ we get
$$\varepsilon (e_{i},e_{j})\varepsilon (e_{j},e_{i})^{-1}=q_{ij}
\ \ \ \mbox{ for }i\ne j.$$   For $1\le i\ne j\le l,\; m,n\in \Z$,
we have
\begin{eqnarray*}
& &\bar{X}^{\varepsilon}_{i,m}\bar{X}^{\varepsilon}_{j,n}
=\varepsilon(e_{i},e_{j})(\bar{X}_{i,m}\bar{X}_{j,n})^{\varepsilon}
=\varepsilon(e_{i},e_{j})\varepsilon(e_{j},e_{i})^{-1}
\bar{X}^{\varepsilon}_{j,n}\bar{X}^{\varepsilon}_{i,m}
=q_{ij}\bar{X}^{\varepsilon}_{j,n}\bar{X}^{\varepsilon}_{i,m},\\
 & &\bar{Y}^{\varepsilon}_{i,m}\bar{Y}^{\varepsilon}_{j,n}
=q_{ij}\bar{Y}^{\varepsilon}_{j,n}\bar{Y}^{\varepsilon}_{i,m},\\
& &\bar{X}^{\varepsilon}_{i,m}\bar{Y}^{\varepsilon}_{j,n}
=\varepsilon(e_{i},-e_{j})\varepsilon(-e_{j},e_{i})^{-1}
\bar{Y}^{\varepsilon}_{j,n}\bar{X}^{\varepsilon}_{i,m}
=q_{ij}^{-1}\bar{Y}^{\varepsilon}_{j,n}\bar{X}^{\varepsilon}_{i,m}
=q_{ji}\bar{Y}^{\varepsilon}_{j,n}\bar{X}^{\varepsilon}_{i,m}.
\end{eqnarray*}
Since $\varepsilon(e_{i},e_{i})=1$, we also have
\begin{eqnarray*}
& &\bar{X}^{\varepsilon}_{i,m}\bar{X}^{\varepsilon}_{i,n}
=(\bar{X}_{i,n}\bar{X}_{i,m})^{\varepsilon}
=q_{ii}(\bar{X}_{i,m}\bar{X}_{i,n})^{\varepsilon}
=q_{ii}\bar{X}^{\varepsilon}_{i,n}\bar{X}^{\varepsilon}_{i,m},\\
& &\bar{Y}^{\varepsilon}_{i,m}\bar{Y}^{\varepsilon}_{i,n}
=q_{ii}\bar{Y}^{\varepsilon}_{i,n}\bar{Y}^{\varepsilon}_{i,m},\\
& &\bar{X}^{\varepsilon}_{i,m}\bar{Y}^{\varepsilon}_{i,n}
-q_{ii}\bar{Y}^{\varepsilon}_{i,n}\bar{X}^{\varepsilon}_{i,m}
=(\bar{X}_{i,m}\bar{Y}_{i,n})^{\varepsilon}
-q_{ii}(\bar{Y}_{i,n}\bar{X}_{i,m})^{\varepsilon} =\delta_{m+n+1,0}.
\end{eqnarray*}
 It follows that there exists an
algebra homomorphism $\psi$ from $\A_{\bf Q}$ onto
$\A^{\varepsilon}$ such that
$$\psi(X_{i,m})=\bar{X}_{i,m}^{\varepsilon},\ \ \ \
\psi(Y_{i,m})=\bar{Y}_{i,m}^{\varepsilon}\ \ \ \ \mbox{ for
}i=1,\dots,l,\; m\in \Z.$$

On the other hand, define a group homomorphism $\varepsilon':
\Z^{l}\times \Z^{l}\rightarrow \C^{\times}$ by
$$\varepsilon'(e_{i},e_{j})=\varepsilon(e_{i},e_{j})^{-1}\ \ \
\mbox{  for }1\le i,j\le l.$$ We then consider the
$\varepsilon'$-twist $\A_{\bf Q}^{\varepsilon'}$ of $\A_{\bf Q}$. It
is straightforward to show that for $1\le i\ne j\le l,\; m,n\in \Z$,
\begin{eqnarray*}
& &X^{\varepsilon'}_{i,m}X^{\varepsilon'}_{j,n}
=X^{\varepsilon'}_{j,n}X^{\varepsilon'}_{i,m},\ \
Y^{\varepsilon'}_{i,m}Y^{\varepsilon'}_{j,n}
=Y^{\varepsilon'}_{j,n}Y^{\varepsilon'}_{i,m},\ \
X^{\varepsilon'}_{i,m}Y^{\varepsilon'}_{j,n}
=Y^{\varepsilon'}_{j,n}X^{\varepsilon'}_{i,m},\\
& &X^{\varepsilon'}_{i,m}X^{\varepsilon'}_{i,n}
=q_{ii}X^{\varepsilon'}_{i,n}X^{\varepsilon'}_{i,m},\ \
Y^{\varepsilon'}_{i,m}Y^{\varepsilon'}_{i,n}
=q_{ii}Y^{\varepsilon'}_{i,n}Y^{\varepsilon'}_{i,m},\ \
X^{\varepsilon'}_{i,m}Y^{\varepsilon'}_{i,n}
-q_{ii}Y^{\varepsilon'}_{i,n}X^{\varepsilon'}_{i,m}
 =\delta_{m+n+1,0}.
\end{eqnarray*}
It follows that there exists an algebra homomorphism from $\A$ to
$\A_{\bf Q}^{\varepsilon'}$, sending $\bar{X}_{i,m}$ to
$X_{i,m}^{\varepsilon'}$ and $\bar{Y}_{i,m}$ to
$Y_{i,m}^{\varepsilon'}$ for $1\le i\le l,\; m\in \Z$. Furthermore,
there exists an algebra homomorphism from $\A^{\varepsilon}$ to
$(\A_{\bf Q}^{\varepsilon'})^{\varepsilon}=\A_{\bf Q}$,
 sending $\bar{X}_{i,m}^{\varepsilon}$ to
$X_{i,m}$ and $\bar{Y}_{i,m}^{\varepsilon}$ to $Y_{i,m}$ for $1\le
i\le l,\; m\in \Z$. Now it follows that $\psi$ is an isomorphism
from $\A_{\bf Q}$ onto $\A^{\varepsilon}$. Notice that for
homogeneous $a_{1}\in \A_{(q_{11})}[m_{1}],\dots,a_{l}\in
\A_{(q_{ll})}[m_{l}]$, we have
$$a_{1}^{\varepsilon}\cdots a_{l}^{\varepsilon}=\lambda (a_{1}\cdots
a_{l})^{\varepsilon}$$ in $\A^{\varepsilon}$ for some nonzero number
$\lambda$ depending on $m_{1},\dots,m_{l}$. Then all the assertions
follow immediately.
\end{proof}

Furthermore we have:

\bc{cpbwforaq} The linear maps
\begin{eqnarray}
&\A_{(q_{11})}^{\pm}\otimes \cdots \otimes
\A_{(q_{ll})}^{\pm} & \rightarrow \A_{\bf Q}^{\pm}\nonumber\\
&a_{1}\otimes \cdots \otimes a_{l}& \mapsto a_{1}\cdots a_{l}
\end{eqnarray}
are linear isomorphisms.  Furthermore, the linear map
$$\A_{\bf Q}^{-}\otimes
\A_{\bf Q}^{+}\rightarrow \A_{\bf Q};\ \ a\otimes b\mapsto ab$$ is a
linear isomorphism.
 \ec

\begin{proof} For each $1\le i\le l$, it is well known that the linear map
$$\A_{(q_{ii})}^{-}\otimes \A_{(q_{ii})}^{+}\rightarrow
\A_{(q_{ii})};\; a\otimes b\mapsto ab$$ is a linear isomorphism.
With $\A=\A_{(q_{11})}\otimes\cdots \otimes \A_{(q_{ll})}$, the
linear map
$$\A^{-}\otimes \A^{+}\rightarrow \A;\ \ a\otimes
b\mapsto ab$$ is a linear isomorphism,  where
$$\A^{\pm}=\A_{(q_{11})}^{\pm}\otimes \cdots \otimes \A_{(q_{ll})}^{\pm}.$$
Clearly, $\A^{\pm}$ are $\Z^{l}$-graded subalgebras. For any
homogeneous vectors $a_{1},\dots, a_{k}\in \A$, we have
$$a_{1}^{\varepsilon}\cdots a_{k}^{\varepsilon}
=\mu(a_{1}\cdots a_{k})^{\varepsilon},$$ where $\mu$ is a nonzero
number depending only on the gradings of $a_{1},\dots,a_{k}$. It
follows that $\A_{\bf Q}^{\pm}=(\A^{\pm})^{\varepsilon}$, under the
identification of $\A_{\bf Q}$ with $\A^{\varepsilon}$ in
Proposition \ref{ptwisted-algebra}. Then all the assertions follow
{}from the second assertion of Proposition \ref{ptwisted-algebra}.
\end{proof}

In the following we shall show that $\A_{\bf Q}$ can also be
constructed inductively by using smash products of algebras
$\A_{(q_{11})},\dots,\A_{(q_{ll})}$.  First we formulate the
following straightforward consequence of the definition:

\bl{lsigma} For any ${\bf q}=(q_{1},\dots,q_{l})\in
(\C^{\times})^{l}$, there exists a (unique) automorphism
$\sigma_{\bf q}$ of $\A_{\bf Q}$ such that
\begin{eqnarray}
\sigma_{\bf q} (X_{i,m})=q_{i}X_{i,m},\ \ \ \ \sigma_{\bf q}
(Y_{i,m})=q_{i}^{-1}Y_{i,m}
\end{eqnarray}
for $i=1,\dots,l,\; m\in \Z$. Furthermore, the multiplicative group
$(\C^{\times})^{l}$ acts on $\A_{\bf Q}$ with ${\bf
q}=(q_{1},\dots,q_{l})\in (\C^{\times})^{l}$ acting as $\sigma_{\bf
q}$. \el

Note that any associative algebra with identity is a nonlocal vertex
algebra. Then all the results for smash product nonlocal vertex
algebras in Section 2 hold for the degenerate case with $U=A$ and
$V=B$ classical associative algebras.

 \bp{piterated} Let ${\bf
 Q}=(q_{ij})_{i,j=1}^{l}$ be a complex matrix with $l\ge 2$ such that
  $q_{ij}q_{ji}=1$ for $1\le i,j\le l$. Set ${\bf Q}'=(q_{ij})_{1\le
  i,j\le l-1}$ and set
$${\bf q}=(q_{1,l},\dots,q_{l-1,l})\in (\C^{\times})^{l-1}.$$
Let $\Z$ act on $\A_{{\bf Q}'}$ with $n$ acting as $\sigma_{\bf
q}^{n}$ for $n\in \Z$ and equip $\A_{(q_{ll})}$ with the
$\Z$-grading constructed before. Then the natural algebra
homomorphisms from $\A_{{\bf Q}'}$ and from $\A_{(q_{ll})}$ to
$\A_{\bf Q}$ give rise to an algebra isomorphism from $\A_{{\bf
Q}'}\sharp_{\Z} \A_{(q_{ll})}$ to $\A_{\bf Q}$. \ep

\begin{proof} Let $X_{m},Y_{m}$ $(m\in \Z)$ denote
the generators of $\A_{(q_{ll})}$.
For $1\le i,j\le l-1,\; m\in \Z$, set
  $\tilde{X}_{i,m}=X_{i,m},\;\tilde{Y}_{i,m}=Y_{i,m}$, and
$\tilde{X}_{l,m}=X_{m},\; \tilde{Y}_{l,m}=Y_{m}$, which are elements
of $\A_{{\bf Q}'}\sharp_{\Z} \A_{(q_{ll})}$. For $1\le i\le l-1,\;
m,n\in \Z$, noticing that $\deg X_{m}=-1$ and $\deg Y_{n}=1$, we
have
\begin{eqnarray*}
& &\tilde{X}_{l,n}\tilde{X}_{i,m}=X_{n}X_{i,m} =\sigma_{\bf
q}^{-1}(X_{i,m})X_{n}=q_{i,l}^{-1}X_{i,m}X_{n}=
q_{l,i}\tilde{X}_{i,m}\tilde{X}_{l,n},\\
& &\tilde{Y}_{l,n}\tilde{Y}_{i,m}=Y_{n}Y_{i,m} =\sigma_{\bf
q}(Y_{i,m})Y_{n}=q_{i,l}^{-1}Y_{i,m}Y_{n}=
q_{l,i}\tilde{Y}_{i,m}\tilde{Y}_{l,n},\\
& &\tilde{Y}_{l,n}\tilde{X}_{i,m}=Y_{n}X_{i,m} =\sigma_{\bf
q}(X_{i,m})Y_{n}=q_{i,l}X_{i,m}Y_{n}
=q_{l,i}^{-1}\tilde{X}_{i,m}\tilde{Y}_{l,n},\\
& &\tilde{X}_{l,n}\tilde{Y}_{i,m}=X_{n}Y_{i,m} =\sigma_{\bf
q}^{-1}(Y_{i,m})X_{n}=q_{i,l}Y_{i,m}X_{n}
=q_{i,l}\tilde{Y}_{i,m}\tilde{X}_{l,n}.
\end{eqnarray*}
It follows that there exists an algebra homomorphism $\pi$ from
$\A_{\bf Q}$ to $\A_{{\bf Q}'}\sharp_{\Z} \A_{(q_{ll})}$, sending
$X_{i,m}$ to $\tilde{X}_{i,m}$, $Y_{i,m}$ to $\tilde{Y}_{i,m}$ for
$1\le i,j\le l,\; m\in \Z$. On the other hand, let $f$ be the
natural algebra homomorphism from $\A_{{\bf Q}'}$ to $\A_{\bf Q}$
and $g$ the natural algebra homomorphism {}from $\A_{(q_{ll})}$ to
$\A_{\bf Q}$. By Proposition \ref{psmash-universal} we have  an
algebra homomorphism $f\sharp g$ from $\A_{{\bf Q}'}\sharp_{\Z}
\A_{(q_{ll})}$ to $\A_{\bf Q}$, extending $f$ and $g$. It is clear
that $f\sharp g$ is inverse to $\pi$, so that $f\sharp g$ is an
isomorphism.
\end{proof}

Recall that $\A_{\bf Q}^{\pm}$ are the subalgebras of $\A_{\bf Q}$,
generated by $X_{i,m},Y_{i,m}$ for $1\le i\le l,\; m\ge 0$. A vector
$w$ in an $\A_{\bf Q}$-module $W$ is called a {\em vacuum vector} if
$\A_{\bf Q}^{+}w=0$, and  an $\A_{\bf Q}$-module $W$ equipped with a
vacuum vector which generates $W$ is called a {\em vacuum $\A_{\bf
Q}$-module}.

Set
\begin{eqnarray}
V_{\bf Q}=\A_{\bf Q}/(\A_{\bf Q}\A_{\bf Q}^{+}),
\end{eqnarray}
a left $\A_{\bf Q}$-module, and set
$${\bf 1}=1+(\A_{\bf Q}\A_{\bf Q}^{+})\in V_{\bf Q}.$$
Clearly, ${\bf 1}$ is a vacuum vector and $V_{\bf Q}$ equipped with
${\bf 1}$
 is a vacuum $\A_{\bf Q}$-module. Furthermore,
$V_{\bf Q}$ is universal in the obvious sense.

\bp{pAQ}  The vacuum $\A_{\bf Q}$-module $V_{\bf Q}$ is irreducible
and  every nonzero vacuum $\A_{\bf Q}$-module is irreducible and
isomorphic to $V_{\bf Q}$. \ep

\begin{proof} Notice that if the
universal vacuum module $V_{\bf Q}$ is irreducible, the other
assertions follow immediately. We shall use induction on $l$ to show
that $V_{\bf Q}$ is an irreducible $\A_{\bf Q}$-module. If $l=1$,
$\A_{\bf Q}=\A_{(q_{11})}$ is a Weyl algebra or a Clifford algebra
and it is well known that
 $V_{(q_{11})}$ is an irreducible $\A_{(q_{11})}$-module.
 Now, assume that $l\ge 2$. For a given $l\times l$ matrix ${\bf Q}$,
 set ${\bf
Q}'=(q_{ij})_{i,j=1}^{l-1}$. By the induction hypothesis $V_{{\bf
Q}'}$ is an irreducible $\A_{{\bf Q}'}$-module and every nonzero
vacuum $\A_{{\bf Q}'}$-module is isomorphic to $V_{{\bf Q}'}$.
Recall that $\A_{\bf Q}=\A_{{\bf Q}'}\sharp \A_{(q_{ll})}$
(Proposition \ref{piterated}).  By Proposition
\ref{pproduct-module}, $V_{{\bf Q}'}\sharp V_{(q_{ll})}$ is an
irreducible $\A_{{\bf Q}'}\sharp \A_{(q_{ll})}$-module. It is clear
that ${\bf 1}\otimes {\bf 1}$ is a vacuum vector, so that $V_{{\bf
Q}'}\sharp V_{(q_{ll})}$ is a vacuum $\A_{\bf Q}$-module.  Let
$\psi$ denote the natural $\A_{\bf Q}$-homomorphism from $\A_{\bf
Q}=\A_{{\bf Q}'}\sharp \A_{(q_{ll})}$ to $V_{\bf Q}$, sending $a\in
\A_{\bf Q}$ to $a{\bf 1}$. Let $a'\in \A_{{\bf Q}'},\; b\in
\A_{(q_{ll})}[m]$ with $m\in \Z$. We have $\psi(a'\sharp b)=a'b{\bf
1}=b\sigma_{\bf q}^{-m}(a'){\bf 1}$. As $\sigma_{\bf q} (\A_{{\bf
Q}'}^{+})\subset \A_{{\bf Q}'}^{+}$, we see that $\psi$ gives rise
to an $\A_{\bf Q}$-homomorphism $\bar{\psi}$ from $V_{{\bf
Q}'}\sharp V_{(q_{ll})}$ to $V_{\bf Q}$, sending ${\bf 1}\otimes
{\bf 1}$ to ${\bf 1}$. As $V_{{\bf Q}'}\sharp V_{(q_{ll})}$ is
irreducible, it follows that $\bar{\psi}$ is an isomorphism.
Consequently, $V_{\bf Q}$ is an irreducible $\A_{\bf Q}$-module.
\end{proof}

For $1\le i\le l$, set
\begin{eqnarray}
u^{(i)}=X_{i,-1}{\bf 1},\ \ v^{(i)}=Y_{i,-1}{\bf 1}\in V_{\bf Q}
\end{eqnarray}
and set
\begin{eqnarray}
X_{i}(x)=\sum_{n\in \Z}X_{i,n}x^{-n-1},\ \ Y_{i}(x)=\sum_{n\in
\Z}Y_{i,n}x^{-n-1}\in \A_{\bf Q}[[x,x^{-1}]].
\end{eqnarray}

\bt{tsecond-main} Let ${\bf Q}=(q_{ij})_{1\le i,j\le l}$ be a
complex matrix such that $q_{ij}q_{ji}=1$ for $1\le i,j\le l$, let
$\A_{\bf Q}$ be the associative algebra associated with ${\bf Q}$
and let $V_{\bf Q}$ be the universal vacuum $\A_{\bf Q}$-module.
There exists a (unique) irreducible quantum vertex algebra structure
on $V_{\bf Q}$ with ${\bf 1}$ as the vacuum vector such that
$$Y(u^{(i)},x)=X_{i}(x),\ \
Y(v^{(i)},x)=Y_{i}(x)\ \ \ \mbox{ for }i=1,\dots,l.$$ Let $W$ be any
$\A_{\bf Q}$-module satisfying the condition that for any $w\in W$,
$X_{i,m}w=Y_{i,m}w=0$ for $1\le i\le l$ and for $m$ sufficiently
large. Then there exists a (unique) $V_{\bf Q}$-module structure on
$W$ with
$$Y_{W}(u^{(i)},x)=X_{i}(x),\ \ \ Y_{W}(v^{(i)},x)=Y_{i}(x)
\ \ \ \mbox{ for }i=1,\dots,l.$$ Conversely, any $V_{\bf Q}$-module
$W$ is an $\A_{\bf Q}$-module with
$$X_{i}(x)=Y_{W}(u^{(i)},x),\ \ \ Y_{i}(x)=Y_{W}(v^{(i)},x)
\ \ \ \mbox{ for }i=1,\dots,l.$$
 \et

\begin{proof} For $m\in \Z,\; k\ge 1$, let
$$Z_{m},Z^{(1)}_{m},\dots,Z^{(k)}_{m}\in \{
X_{i,m},Y_{i,m}\;|\; 1\le i\le l\}.$$ It follows from induction on
$k$ that
$$Z_{m}Z^{(1)}_{m_{1}}\cdots Z^{(k)}_{m_{k}}{\bf 1}=0$$
whenever $m\ge 0$ and $m\notin \{-m_{1}-1,\dots,-m_{k}-1\}$. Then
$$\{X_{i}(x),Y_{i}(x)\;|\;1\le i\le l\}\subset
\E(V_{\bf Q})=\Hom (V_{\bf Q},V_{\bf
Q}((x))).$$ Set $T=\{ X_{i}(x),Y_{i}(x)\;|\; i=1,\dots,l\}$. {}From
the defining relations of $\A_{\bf Q}$, $T$ is an $\S$-local subset
of $\E(V_{\bf Q})$. By Theorem 5.8 of \cite{li-qva1}, $T$ generates
a weak quantum vertex algebra
 $\<T\>$, where the vertex operator map is denoted by $Y_{\E}$,
 with $V_{\bf Q}$ as a module where
 $Y_{V_{\bf Q}}(\alpha(x),x_{0})=\alpha(x_{0})$ for $\alpha(x)\in \<T\>$.
By Proposition 6.7 of \cite{li-qva1}, we have
\begin{eqnarray*}
& &Y_{\E}(X_{i}(x),x_{1})Y_{\E}(X_{j}(x),x_{2})
=q_{ij}Y_{\E}(X_{j}(x),x_{2})Y_{\E}(X_{i}(x),x_{1}),\\
& &Y_{\E}(Y_{i}(x),x_{1})Y_{\E}(Y_{j}(x),x_{2})
=q_{ij}Y_{\E}(Y_{j}(x),x_{2})Y_{\E}(Y_{i}(x),x_{1}),\\
& &Y_{\E}(X_{i}(x),x_{1})Y_{\E}(Y_{j}(x),x_{2})
-q_{ji}Y_{\E}(Y_{j}(x),x_{2})Y_{\E}(X_{i}(x),x_{1})
=\delta_{ij}x_{2}^{-1}\delta\left(\frac{x_{1}}{x_{2}}\right)
\end{eqnarray*}
for $i,j=1,\dots,l$. It follows that $\<T\>$ is an $\A_{\bf
Q}$-module with
$$X_{i}(x_{0})=Y_{\E}(X_{i}(x),x_{0}),\ \
Y_{i}(x_{0})=Y_{\E}(Y_{i}(x),x_{0})\ \ \ \mbox{ for }1\le i\le l.$$
It is clear that $\<T\>$ is a vacuum module with $1_{V_{\bf Q}}$ as
the vacuum vector. From the construction of $V_{\bf Q}$, there is an
$\A_{\bf Q}$-module homomorphism $\pi$ from $V_{\bf Q}$ to $\<T\>$,
sending ${\bf 1}$ to $1_{V_{\bf Q}}$. Now, in view of Theorem 6.3 of
\cite{li-qva1}, there exists a weak quantum vertex algebra structure
on $V_{\bf Q}$ with all the required properties. As $V_{\bf Q}$ is
an irreducible $\A_{\bf Q}$-module, $V_{\bf Q}$ as a $V_{\bf
Q}$-module is irreducible. Therefore, $V_{\bf Q}$ is an irreducible
quantum vertex algebra.

Let $W$ be any $\A_{\bf Q}$-module satisfying the condition that for
any $w\in W$, $X_{i,m}w=Y_{i,m}w=0$ for $1\le i\le l$ and for $m$
sufficiently large. Then $T_{W}=\{ X_{i}(x),Y_{i}(x)\;|\;
i=1,\dots,l\}$, viewed as a subset of $\E(W)$,
 is $\S$-local  and then by Theorem 5.8 of \cite{li-qva1},
$T_{W}$ generates a weak quantum vertex algebra $\<T_{W}\>$
 with $W$ as a module. It follows from the same argument that
$\<T_{W}\>$ is a vacuum $\A_{\bf Q}$-module and that there exists an
$\A_{\bf Q}$-module homomorphism $\psi$ from $V_{\bf Q}$ to
$\<T_{W}\>$, sending ${\bf 1}$ to $1_{W}$. As $V_{\bf Q}$ is
generated by $u^{(i)},v^{(i)}$, it follows that $\psi$ is a
weak-quantum-vertex-algebra homomorphism. Consequently, $W$ is a
$V_{\bf Q}$-module.

On the other hand, let $(W,Y_{W})$ be a $V_{\bf Q}$-module. By
Proposition 6.7 of \cite{li-qva1}, we have
\begin{eqnarray*}
& &Y_{W}(u^{(i)},x_{1})Y_{W}(u^{(j)},x_{2})
=q_{ij}Y_{W}(u^{(j)},x_{2})Y_{W}(u^{(i)},x_{1}),\\
& &Y_{W}(v^{(i)},x_{1})Y_{W}(v^{(j)},x_{2})
=q_{ij}Y_{W}(v^{(j)},x_{2})Y_{W}(v^{(i)},x_{1}),\\
& &Y_{W}(u^{(i)},x_{1})Y_{W}(v^{(j)},x_{2})
-q_{ji}Y_{W}(v^{(j)},x_{2})Y_{W}(u^{(i)},x_{1})
=\delta_{ij}x_{2}^{-1}\delta\left(\frac{x_{1}}{x_{2}}\right)
\end{eqnarray*}
for $1\le i,j\le l$. It follows that $W$ is an $\A_{\bf Q}$-module
with $X_{i}(x)=Y_{W}(u^{(i)},x)$ and $Y_{i}(x)=Y_{W}(v^{(i)},x)$ for
$i=1,\dots,l$. For any $w\in W$, as $Y_{W}(a,x)w\in W((x))$ for
$a\in V_{\bf Q}$, we have $X_{i,m}w=Y_{i,m}w=0$ for $1\le i\le l$
and for $m$ sufficiently large.
\end{proof}

Furthermore, we have:

\bp{puniversal-vq} Let $V$ be any nonlocal vertex algebra and let
$\psi$ be any map {}from $\{ u^{(i)},v^{(i)}\;|\;i=1,\dots,l\}$ to
$V$ such that
\begin{eqnarray*}
& &Y(\psi(u^{(i)}),x_{1})Y(\psi(u^{(j)}),x_{2}) =
q_{ij}Y(\psi(u^{(j)}),x_{2})Y(\psi(u^{(i)}),x_{1}),\\
& &Y(\psi(v^{(i)}),x_{1})Y(\psi(v^{(j)}),x_{2}) = q_{ij}Y(
\psi(v^{(j)}),x_{2})Y(\psi(v^{(i)}),x_{1}),\\
& &Y(\psi(u^{(i)}),x_{1})Y(\psi(v^{(j)}),x_{2})-
q_{ji}Y(\psi(v^{(j)}),x_{2})Y(\psi(u^{(i)}),x_{1})
=\delta_{ij}x_{2}^{-1}\delta\left(\frac{x_{1}}{x_{2}}\right)\ \ \ \
\end{eqnarray*}
for $1\le i,j\le l$. Then there exists a unique
nonlocal-vertex-algebra homomorphism from $V_{\bf Q}$ to $V$,
extending $\psi$. \ep

\begin{proof} The uniqueness is clear.
{}From the assumption we see that $V$ is an $\A_{\bf Q}$-module with
$X_{i}(x)=Y(\psi(u^{(i)}),x)$ and $Y_{i}(x)=Y(\psi(v^{(i)}),x)$ for
$1\le i\le l$, that the submodule generated by the vacuum vector is
a vacuum $\A_{\bf Q}$-module. Then there exists an $\A_{\bf
Q}$-module homomorphism $\tilde{\psi}$ {}from $V_{\bf Q}$ to $V$,
sending the vacuum vector to the vacuum vector. As $V_{\bf Q}$ is
generated by $ u^{(i)},v^{(i)}$ ($i=1,\dots,l$), it follows that
$\tilde{\psi}$ is a nonlocal-vertex-algebra homomorphism, extending
$\psi$.
\end{proof}

Let $V$ be a nonlocal vertex algebra. Define $LC(V)$ to consist of
vectors $a\in V$ satisfying the condition that for every $v\in V$,
there exists a nonnegative integer $k$ such that
$$(x_{1}-x_{2})^{k}[Y(a,x_{1}),Y(v,x_{2})]=0.$$
One can show that $LC(V)$ is a subalgebra of $V$ and it is a vertex
algebra itself.

 \bd{dcore} {\em Let $V$ be a nonlocal vertex algebra.
A vector $\omega\in V$ is called a {\em conformal vector} if
$\omega\in LC(V)$ and if the following Virasoro algebra relation
holds for $m,n\in \Z$,
\begin{eqnarray*}
[L(m),L(n)]=(m-n)L(m+n)+\frac{1}{12}(m^{3}-m)\delta_{m+n,0}c
\end{eqnarray*}
and
\begin{eqnarray} [L(-1),Y(v,x)]=Y(L(-1)v,x)=\frac{d}{dx}Y(v,x)
\end{eqnarray}
for all $v\in V$, where $Y(\omega,x)=\sum_{n\in \Z}L(n)x^{-n-2}$ and
$c$ is a complex number. A nonlocal vertex algebra equipped with a
conformal vector is called a {\em conformal nonlocal vertex
algebra}. } \ed

 \br{rrecall} {\em Recall that for $\epsilon=\pm 1$,
$\A_{(\epsilon)}$ is the associative algebra with generators
$a_{m},b_{n}$ for $m,n\in \Z$, subject to relations
$$a_{m}a_{n}=\epsilon
a_{n}a_{m},\ \ \ b_{m}b_{n}=\epsilon b_{n}b_{m},\ \ \
a_{m}b_{n}-\epsilon b_{n}a_{m}=\delta_{m+n+1,0}.$$  Set
$a=a_{-1}{\bf 1},\; b=b_{-1}{\bf 1}\in V_{(\epsilon)}$. It has been
known (\cite{ffr},  \cite{wei}; cf. \cite{duncan1}, \cite{duncan2})
that there exists a (unique) conformal vertex superalgebra structure
on $V_{(\epsilon)}$ with ${\bf 1}$ as the vacuum vector and with
$$Y(a,x)=\sum_{n\in \Z}a_{n}x^{-n-1}\ \ \mbox{ and }\ \ \
Y(b,x)=\sum_{n\in \Z}b_{n}x^{-n-1},$$ and with  the conformal vector
\begin{eqnarray} \omega=\frac{1}{2}(b_{-2}a-\epsilon a_{-2}b)
\end{eqnarray}
of central charge $-\epsilon$.  Furthermore, $V_{(\epsilon)}$ is
simple and is $\frac{1}{2}\Z_{+}$-graded by the $L(0)$-weights with
$$(V_{(\epsilon)})_{(0)}=\C {\bf 1}\ \ \mbox{ and }\ \
(V_{(\epsilon)})_{(1/2)}=\C a+\C b.$$ On the other hand,
$\A_{(\epsilon)}$ is a $\Z$-graded algebra with $\deg a_{m}=1,\ \deg
b_{n}=-1$ for $m,n\in \Z$, and the vacuum module $V_{(\epsilon)}$ is
naturally a $\Z$-graded $\A_{(\epsilon)}$-module with $\deg {\bf
1}=0$. For $m,n\in \Z$, we have
$$a_{m}V_{(\epsilon)}[n]\subset V_{(\epsilon)}[n+1],\ \ \ \
b_{m}V_{(\epsilon)}[n]\subset V_{(\epsilon)}[n-1].$$  Since
$V_{(\epsilon)}$ as a vertex superalgebra is generated by $a$ and
$b$, it follows from induction that $V_{(\epsilon)}$ equipped with
the $\Z$-grading is a vertex $\Z$-graded superalgebra with
$\omega\in (V_{(\epsilon)})[0]$.} \er

\bp{pvirasoro} Let ${\bf Q}$ be an $l\times l$ complex matrix as
before. Then $V_{\bf Q}$ is a conformal quantum vertex algebra of
central charge $-(q_{11}+\cdots +q_{ll})$ with conformal vector
\begin{eqnarray}
\omega=\frac{1}{2}\sum_{i=1}^{l}(v^{(i)}_{-2}u^{(i)}
-q_{ii}u^{(i)}_{-2}v^{(i)})
\end{eqnarray}
and $V_{\bf Q}$ is $\frac{1}{2}\Z_{+}$-graded by the $L(0)$-weights
such that $(V_{\bf Q})_{(0)}=\C{\bf 1}$ and
$$(V_{\bf Q})_{(1/2)}={\rm span}\{u^{(i)},v^{(i)}\;|\; 1\le
i\le l\}.$$ \ep

\begin{proof}
We have simple conformal vertex superalgebras
$V_{(q_{11})},\dots,V_{(q_{ll})}$.  For $1\le i\le l$, denote the
two generators of $V_{(q_{ii})}$ by $a^{(i)},b^{(i)}$. Set
$$U=V_{(q_{11})}\otimes \cdots \otimes V_{(q_{ll})},$$
the tensor product conformal vertex superalgebra with the conformal
vector
$$\omega=\frac{1}{2}\sum_{i=1}^{l}(b^{(i)}_{-2}a^{(i)}
-q_{ii}a^{(i)}_{-2}b^{(i)}).$$ Then $U$ is
$\frac{1}{2}\Z_{+}$-graded by the $L(0)$-weights with
$U_{(0)}=\C{\bf 1}$ and
$$U_{(1/2)}={\rm span}\{ a^{(i)},b^{(i)}\;|\; 1\le i\le l\}.$$
On the other hand, as $V_{(q_{11})},\dots,V_{(q_{ll})}$ are vertex
$\Z$-graded superalgebras, $U$ is a vertex $\Z^{l}$-graded
superalgebra.
 Let $\varepsilon:
\Z^{l}\times \Z^{l}\rightarrow \C^{\times}$ be the group
homomorphism defined in (\ref{def-varepsilon}). In view of Lemma
\ref{ltwisted}, we have a nondegenerate quantum vertex algebra
$(U^{\varepsilon},Y_{\varepsilon},{\bf 1})$. A straightforward
calculation yields
\begin{eqnarray*}
& &Y_{\varepsilon}(a^{(i)},x_{1})Y_{\varepsilon}(a^{(j)},x_{2}) =
q_{ij}Y_{\varepsilon}(a^{(j)},x_{2})Y_{\varepsilon}(a^{(i)},x_{1}),\\
& &Y_{\varepsilon}(b^{(i)},x_{1})Y_{\varepsilon}(b^{(j)},x_{2}) =
q_{ij}Y_{\varepsilon}(b^{(j)},x_{2})Y_{\varepsilon}(b^{(i)},x_{1}),\\
& &Y_{\varepsilon}(a^{(i)},x_{1})Y_{\varepsilon}(b^{(j)},x_{2})-
q_{ji}Y_{\varepsilon}(b^{(j)},x_{2})Y_{\varepsilon}(a^{(i)},x_{1})
=\delta_{ij}x_{2}^{-1}\delta\left(\frac{x_{1}}{x_{2}}\right)\ \ \ \
\end{eqnarray*}
for $1\le i,j\le l$. Hence, by Proposition \ref{puniversal-vq},
there is a nonlocal vertex algebra homomorphism $\psi$ from $V_{\bf
Q}$ to $U^{\varepsilon}$ such that
 $$\psi(u_{i})=(a^{(i)})^{\varepsilon},\ \ \ \
 \psi(v^{i})=(b^{(i)})^{\varepsilon}\ \ \ \mbox{ for }1\le i\le l.$$
 As $V_{\bf Q}$ is simple, $\psi$ is an isomorphism.
 Notice that $\omega\in U[\bf 0]$, the degree-${\bf 0}$ subspace of
 $U$ with respect to the $\Z^{l}$-grading. It follows that $\omega\in
 LC(U^{\varepsilon})$ and
 $Y_{\varepsilon}(\omega,x)=Y(\omega,x)$ on $U$. Then
 $U^{\varepsilon}$ equipped with the conformal vector
 $\omega$ is a conformal vertex algebra.
 \end{proof}

\bp{pvq-module} Let $W$ be a $V_{\bf Q}$-module satisfying the
condition that for any $w\in W$, there exists a positive integer $k$
such that $(\A_{\bf Q}^{+})^{k}w=0$. Then $W$ is a direct sum of
irreducible modules isomorphic to the adjoint module $V_{\bf Q}$.
\ep

\begin{proof} If $l=1$,  $\A_{\bf Q}=\A_{(q_{11})}$
is either a Weyl algebra or a Clifford algebra and it is well known
that the assertion is true. If $l\ge 2$, by Proposition
\ref{piterated}, $\A_{\bf Q}=\A_{{\bf Q}'}\sharp \A_{(q_{ll})}$ with
${\bf Q}'=(q_{ij})_{i,j=1}^{l-1}$. Assume that the assertion holds
for $\A_{{\bf Q}'}$. Let $W$ be any $\A_{\bf Q}$-module satisfying
the condition that for any $w\in W$, there exists a positive integer
$k$ such that $(\A_{\bf Q}^{+})^{k}w=0$. Then $W$ as an $\A_{{\bf
Q}'}$-module is a direct sum of irreducible $\A_{{\bf
Q}'}$-submodules isomorphic to $V_{{\bf Q}'}$. In view of
Proposition \ref{pmodule-class}, we have $W=V_{{\bf Q}'}\sharp F$
for some $\A_{(q_{ll})}$-module $F$. We know that $F$ is a direct
sum of irreducible $\A_{(q_{ll})}$-modules isomorphic to
$V_{(q_{ll})}$. By Proposition \ref{pproduct-module}, $V_{{\bf
Q}'}\sharp V_{(q_{ll})}$ is an irreducible $\A_{{\bf Q}'}\sharp
\A_{(q_{ll})}$-module. It follows that $V_{\bf Q}\simeq V_{{\bf
Q}'}\sharp V_{(q_{ll})}$. Consequently, $W$ is a direct sum of
irreducible modules isomorphic to $V_{\bf Q}$.
\end{proof}

For the rest of this section, we establish certain results for
$V_{\bf Q}$, which we shall need in the next section.

\bd{dpseudo-auto} {\em Let $V$ be a nonlocal vertex algebra. A {\em
pseudo-endomorphism} of $V$ (see \cite{li-pseudo}) is a linear map
$\Phi(x): V\rightarrow V\otimes \C((x))$ such that $\Phi(x)({\bf
1})={\bf 1}$ and
\begin{eqnarray}
\Phi(x_{1})Y(v,x_{2})=Y(\Phi(x_{1}-x_{2})v,x_{2})\Phi(x_{1}) \ \ \ \
\mbox{ for all }v\in V.
\end{eqnarray}
A pseudo-endomorphism $\Phi(x)$ is called a {\em
pseudo-automorphism} if there exists a pseudo-endomorphism $\Psi(x)$
such that $\Phi(x)\Psi(x)v=v=\Psi(x)\Phi(x)v$ for $v\in V$. We say
that pseudo-endomorphisms $\Phi(x)$ and $\Psi(x)$ {\em commute} if
$\Phi(x_{1})\Psi(x_{2})=\Psi(x_{2})\Phi(x_{1})$.}\ed

The following is straightforward (see \cite{li-pseudo}):

 \bp{ppseodu-auto} Let
$V$ be a nonlocal vertex algebra. A pseudo-endomorphism of $V$
exactly amounts to a nonlocal-vertex-algebra homomorphism from $V$
to the tensor product nonlocal vertex algebra $V\otimes \C((x))$,
where $\C((x))$ is viewed as a vertex algebra with $1$ as the vacuum
vector and
$$Y(f(x),z)g(x)=f(x-z)g(x)\ \ \ \mbox{ for }f(x),g(x)\in
\C((x)).$$ \ep

We have the following results for the quantum vertex algebra $V_{\bf
Q}$:

\bl{lpseudo-auto} Let ${\bf Q}=(q_{ij})$ be an $l\times l$ matrix as
before. For any $(p_{1}(x),\dots,p_{l}(x))\in \C((x))^{l}$ with
$p_{i}(x)\ne 0$ for all $i$, there exists a pseudo-automorphism
$\Phi(x)$ of $V_{{\bf Q}}$ such that
\begin{eqnarray*}
\Phi(x)(u^{(i)})=u^{(i)}\otimes p_{i}(x),\ \ \ \ \ \
\Phi(x)(v^{(i)})=v^{(i)}\otimes p_{i}(x)^{-1}\ \ \ \mbox{ for
}i=1,\dots,l.
\end{eqnarray*}
Furthermore, all such pseudo-automorphisms commute.
 \el

\begin{proof} Denote by $\hat{Y}$ the vertex operator map for the
tensor product nonlocal vertex algebra  $V_{\bf Q}\otimes \C((x))$.
{}From definition we have
$$\hat{Y}(a\otimes f(x),z)=Y(a,z)\otimes f(x-z) \ \ \ \ \mbox{ for
}a\in V_{\bf Q},\; f(x)\in \C((x)).$$ The following relations hold:
\begin{eqnarray*}
& &\hat{Y}(u^{(i)}\otimes p_{i}(x),x_{1})\hat{Y}(u^{(j)}\otimes
p_{j}(x),x_{2})\\
&=&Y(u^{(i)},x_{1})Y(u^{(j)},x_{2})\otimes p_{i}(x-x_{1})p_{j}(x-x_{2})\\
&=&q_{ij}Y(u^{(j)},x_{2})Y(u^{(i)},x_{1})
\otimes p_{j}(x-x_{2})p_{i}(x-x_{1})\\
&=&q_{ij}\hat{Y}(u^{(j)}\otimes
p_{j}(x),x_{2})\hat{Y}(u^{(i)}\otimes p_{i}(x),x_{1}),
\end{eqnarray*}
\begin{eqnarray*}
& &\hat{Y}(v^{(i)}\otimes p_{i}(x)^{-1},x_{1})
\hat{Y}(v^{(j)}\otimes p_{j}(x)^{-1},x_{2})\\
&=&q_{ij}\hat{Y}(v^{(j)}\otimes
p_{j}(x)^{-1},x_{2})\hat{Y}(v^{(i)}\otimes p_{i}(x)^{-1},x_{1}),
\end{eqnarray*}
\begin{eqnarray*}
& &\hat{Y}(u^{(i)}\otimes p_{i}(x),x_{1})\hat{Y}(v^{(j)}\otimes
p_{j}(x)^{-1},x_{2})\\
& &\ \ \ \ -q_{ji}\hat{Y}(v^{(j)}\otimes
p_{j}(x)^{-1},x_{2})\hat{Y}(u^{(i)}\otimes p_{i}(x),x_{1})\\
&=&\left(Y(u^{(i)},x_{1})Y(v^{(j)},x_{2})
-q_{ji}Y(v^{(j)},x_{2})Y(v^{(i)},x_{1})\right)
\otimes p_{i}(x-x_{1})p_{j}^{-1}(x-x_{2})\\
&=&\delta_{ij}p_{i}(x-x_{1})p_{j}^{-1}(x-x_{2})
x_{2}^{-1}\delta\left(\frac{x_{1}}{x_{2}}\right)\\
&=&\delta_{ij}p_{i}(x-x_{2})p_{j}^{-1}(x-x_{2})
x_{2}^{-1}\delta\left(\frac{x_{1}}{x_{2}}\right)\\
&=&\delta_{ij}x_{2}^{-1}\delta\left(\frac{x_{1}}{x_{2}}\right).
\end{eqnarray*}
In view of Proposition \ref{puniversal-vq}, there exists a
nonlocal-vertex-algebra homomorphism $\psi$ from $V_{\bf Q}$ into
$V_{\bf Q}\otimes \C((x))$ such that
$$\psi(u^{(i)})=u^{(i)}\otimes p_{i}(x),\ \ \ \
\psi(v^{(i)})=v^{(i)}\otimes p_{i}(x)^{-1} \ \ \ \ \mbox{ for
}i=1,\dots,l.$$ The map $\psi: V_{\bf Q}\rightarrow V_{\bf Q}\otimes
\C((x))$, alternatively denoted by $\Phi(x)$, is a
pseudo-automorphism of $V$ satisfying the required property (recall
Proposition \ref{ppseodu-auto}). The rest is clear.
\end{proof}

\section{Zamolodchikov-Faddeev type quantum vertex algebras}

The associative algebras $\A_{\bf Q}$, which were studied in the
previous section, are actually the simplest Zamolodchikov-Faddeev
algebras. A general Zamolodchikov-Faddeev algebra is associated with
a quantum Yang-Baxter operator $\S(x_{1},x_{2})$ (with two spectral
parameters) on a finite-dimensional vector space. In \cite{li-qva2},
for any finite-dimensional vector space $H$ equipped with a bilinear
form and a linear map $\S(x): H\otimes H\otimes \C[[x]]$, we
constructed a ``universal'' weak quantum vertex algebra $V(H,\S)$
with certain generators and defining relations. Then we studied a
special family of $V(H,\S)$ and proved that $V(H,\S)$ are
nondegenerate quantum vertex algebras. In this section, we study the
case in which $\S(x)$ is diagonalizable. By using the results of
Sections 2 and 3 we prove that either $V(H,S)$ is zero, or an
irreducible quantum vertex algebra with a normal basis in a certain
sense. Furthermore, we prove that $V(H,S)$ are indeed nonzero for a
certain family.

First we recall from \cite{li-qva2} the following notion:

\bd{dQLie-module}
{\em Let $H$ be a vector space
equipped with a bilinear form $\<\cdot,\cdot\>$
and let $\S (x): H\otimes H\rightarrow H\otimes H\otimes \C((x))$ be a
linear map. An {\em $(H,\S)$-module} is
a module $W$ for the (free) tensor algebra $T(H\otimes \C[t,t^{-1}])$
such that
for any $a\in H,\; w\in W$,
\begin{eqnarray}\label{e-truncation-qLie}
a(m)w=0\;\;\;\mbox{ for $m$ sufficiently large},
\end{eqnarray}
where $a(m)$ denotes the operator on $W$ corresponding to $a\otimes
t^{m}$, and such that
\begin{eqnarray}\label{emain-relation-qLie}
a(x_{1})b(x_{2})w
-\sum_{i=1}^{r}f_{i}(x_{2}-x_{1})
b^{(i)}(x_{2})a^{(i)}(x_{1})w
=x_{2}^{-1}\delta\left(\frac{x_{1}}{x_{2}}\right)\<a,b\>w
\end{eqnarray}
for $a,b\in H,\; w\in W$, where $u(x)=\sum_{m\in \Z}u(m)x^{-m-1}$
for $u\in H$ and
$$\S(x)(b\otimes a)=\sum_{i=1}^{r}b^{(i)}\otimes a^{(i)}\otimes
f_{i}(x).$$ (Notice that the condition (\ref{e-truncation-qLie})
guarantees that for any $m,n\in \Z$, the coefficient of
$x_{1}^{m}x_{2}^{n}$ in the left side of (\ref{emain-relation-qLie})
is a finite sum in $W$.)} \ed

Let $W$ be an $(H,\S)$-module. A vector $e\in W$ is called a {\em
vacuum vector}  if
\begin{eqnarray}
a(m)e=0\;\;\;\mbox{ for all }a\in H,\; m\ge 0.
\end{eqnarray}
A {\em vacuum $(H,\S)$-module} is an $(H,\S)$-module $W$ equipped
with a vacuum vector $e$ that generates $W$. We also denote the
vacuum module by $(W,e)$. The notion of a universal vacuum
$(H,\S)$-module is defined in the obvious way. If $\S(x): H\otimes
H\rightarrow H\otimes H\otimes \C[[x]]$, universal vacuum
$(H,\S)$-modules exist,  as was shown in \cite{li-qva2} by the
following tautological construction:
 Let $T(H\otimes
\C[t,t^{-1}])^{+}$ be the subspace of $T(H\otimes \C[t,t^{-1}])$,
spanned by the vectors
$$(a^{(1)}\otimes t^{n_{1}})\cdots (a^{(r)}\otimes t^{n_{r}})$$
for $r\ge 1,\; a^{(i)}\in H,\; n_{i}\in \Z$ with $n_{1}+\cdots +n_{r}\ge 0$.
Set
$$J=T(H\otimes \C[t,t^{-1}])T(H\otimes \C[t,t^{-1}])^{+},$$
a left ideal of $T(H\otimes \C[t,t^{-1}])$, and then set
$$\tilde{V}(H,\S)=T(H\otimes \C[t,t^{-1}])/J,$$
a left $T(H\otimes \C[t,t^{-1}])$-module. Clearly, the condition
(\ref{e-truncation-qLie}) holds. Define $V(H,\S)$ to be the quotient
module of $\tilde{V}(H,\S)$ modulo the relations
(\ref{emain-relation-qLie}). Let ${\bf 1}$ denote the image in
$V(H,\S)$ of $1$.

The following result was established in \cite{li-qva2} (Propositions
2.18 and 4.3):

\bp{pqLie-wqva} Let $H$ be a vector space equipped with a bilinear
form $\<\cdot,\cdot\>$ and let $\S(x): H\otimes H\rightarrow
H\otimes H\otimes \C[[x]]$ be a linear map. Then $(V(H,\S),{\bf 1})$
is a universal vacuum $(H,\S)$-module and there exists a unique weak
quantum vertex algebra structure on $V(H,\S)$ with ${\bf 1}$ as the
vacuum vector such that
\begin{eqnarray*}
Y(a(-1){\bf 1},x)=a(x)\;\;\;\mbox{ for }a\in H.
\end{eqnarray*}
Furthermore, on any $(H,\S)$-module $W$ there exists a unique $V(H,\S)$-module
structure $Y_{W}$ such that
\begin{eqnarray*}
Y_{W}(a(-1){\bf 1},x)=a(x)\;\;\;\mbox{ for }a\in H.
\end{eqnarray*}
\ep

\br{rnot-need} {\em Note that it was assumed in (\cite{li-qva2},
Proposition 2.18) that the linear map from $H$ to $V(H,\S)$, sending
$a$ to $a(-1){\bf 1}$ for $a\in H$, is injective. In fact, this
assumption was not used in the proof. If the linear map $a\in
H\mapsto a(-1){\bf 1}\in V(H,\S)$ is injective, $H$ can be
identified with a subspace of $V(H,\S)$ and $H$ is a generating
subspace. In general, the subspace $\{ a(-1){\bf 1}\;|\; a\in H\}$
is a generating subspace.} \er

\br{runiversal-Vqx} {\em  The nonlocal vertex algebra $V(H,\S)$ is
universal in the sense that for any given nonlocal vertex algebra
$V$ and for any given linear map $\phi: H\rightarrow V$ satisfying
\begin{eqnarray}\label{eyphi}
& &Y(\phi(u),x_{1})Y(\phi(v),x_{2})
-\sum_{i=1}^{r}f_{i}(x_{2}-x_{1})Y(\phi(v^{(i)}),x_{2})Y(\phi(u^{(i)}),x_{1})
\nonumber\\
& &\ \ \ \ \ =\<u,v\>
x_{2}^{-1}\delta\left(\frac{x_{1}}{x_{2}}\right)
\end{eqnarray}
for $u,v\in H$, where
$$\S(x)(v\otimes
u)=\sum_{i=1}^{r}v^{(i)}\otimes u^{(i)}\otimes f_{i}(x)\in H\otimes
H\otimes \C[[x]],$$ there exists a unique nonlocal-vertex-algebra
homomorphism from $V(H,\S)$ to $V$, sending $a(-1){\bf 1}$ to
$\phi(a)$ for $a\in H$. Indeed, the relation (\ref{eyphi}) implies
that $V$ is an $(H,\S)$-module with $a(x)=Y(\phi(a),x)$ for $a\in
H$. Clearly, ${\bf 1}$ is a vacuum vector of $V$ viewed as an
$(H,\S)$-module. It follows that the $(H,\S)$-submodule of $V$
generated from ${\bf 1}$ is a vacuum module. Then there exists an
$(H,\S)$-module homomorphism $\psi$ from $V(H,\S)$ to $V$, sending
the vacuum vector to the vacuum vector. Since $\{ a(-1){\bf
1}\;|\;a\in H\}$ generates $V(H,\S)$ as a nonlocal vertex algebra,
it follows that $\psi$ is a nonlocal-vertex-algebra homomorphism,
where $$\psi(a(-1){\bf 1})=\Res_{x}x^{-1}\psi(Y(a(-1){\bf 1},x){\bf
1})=\Res_{x}x^{-1}Y(\phi(a),x){\bf 1}=\phi(a)$$ for $a\in H$. The
uniqueness is clear.} \er

The following is our first result of this section:

 \bt{tgeneral-zf}
Let $H$ be a finite-dimensional vector space equipped with a
bilinear form $\<\cdot,\cdot\>$ and let $\S(x): H\otimes
H\rightarrow H\otimes H\otimes \C[[x]]$ be a linear map. Suppose
that there exists a basis $\{ u^{(1)},\dots,u^{(l)},
v^{(1)},\dots,v^{(l)}\}$ of $H$ such that
\begin{eqnarray*}
\<u^{(i)},u^{(j)}\>=0=\<v^{(i)},v^{(j)}\>,\ \
\<u^{(i)},v^{(j)}\>=\delta_{ij}=-q_{ii}\<v^{(i)},u^{(j)}\>
\end{eqnarray*}
 and
\begin{eqnarray*}
& &\S(0)(u^{(i)}\otimes u^{(j)})=q_{ji}(u^{(i)}\otimes u^{(j)}),\ \
\ \ \S(0)(v^{(i)}\otimes v^{(j)})=q_{ji}(v^{(i)}\otimes v^{(j)}),\\
& &\S(0)(u^{(i)}\otimes v^{(j)})=q_{ij}(u^{(i)}\otimes v^{(j)}),\ \
\ \ \S(0)(v^{(j)}\otimes u^{(i)})=q_{ij}(v^{(j)}\otimes u^{(i)})
\end{eqnarray*}
where  $q_{ij}\in \C$ with $q_{ij}q_{ji}=1$ for $i,j=1,\dots,l$.
Assume that $V(H,\S)$ is not zero. Then $V(H,\S)$ is an irreducible
quantum vertex algebra with ${\rm Gr}_{F}(V(H,\S))=V_{\bf Q}$ where
${\bf Q}=(q_{ij})_{i,j=1}^{l}$ and $F$ is the increasing filtration
of $V$, defined in Lemma \ref{linduction1}, associated with the
generating subspace $T=\{ a(-1){\bf 1}\;|\; a\in H\}$. Furthermore,
$\S(x)$ must be a unitary rational quantum Yang-Baxter operator. \et

\begin{proof}  We first prove ${\rm Gr}_{F}(V(H,\S))=V_{\bf
Q}$. Recall that $V(H,\S)$ is generated by $a(-1){\bf 1}$ for $a\in
H$, where $Y(a(-1){\bf 1},x)=a(x)$. By definition, for $n<0$,
$F_{n}$ is linearly spanned by the vectors
$$a^{(1)}(m_{1})\cdots a^{(r)}(m_{r}){\bf 1}$$
for $r\ge 1,\; a^{(i)}\in H,\; m_{i}\in \Z$ with $m_{1}+\cdots
+m_{r}\ge -n$, and for $n\ge 0$, $F_{n}$ is linearly spanned by
${\bf 1}$ and the vectors of the same form. By Corollary 4.2 of
\cite{li-qva2}, we have
$$h^{(1)}(m_{1})\cdots h^{(r)}(m_{r}){\bf 1}=0$$
for $r\ge 1,\; h^{(1)},\dots,h^{(r)}\in H,\; m_{1},\dots,m_{r}\in
\Z$ with $m_{1}+\cdots +m_{r}\ge 0$. Consequently, $F_{n}=0$ for
$n\le -1$ and $F_{0}=\C {\bf 1}$ $(\ne 0)$. Then ${\rm
Gr}_{F}(V(H,\S))\ne 0$.

 Notice that $a(-1){\bf 1}\in F_{1}$ for $a\in H$. For $1\le i\le l$, set
$$\bar{u}^{(i)}=u^{(i)}(-1){\bf 1}+F_{0},\ \ \
\bar{v}^{(i)}=v^{(i)}(-1){\bf 1}+F_{0}\in F_{1}/F_{0}\subset {\rm
Gr}_{F}(V(H,\S)).$$ Then ${\rm Gr}_{F}(V(H,\S))$ is generated by
$\bar{u}^{(i)},\bar{v}^{(i)}$ for $1\le i\le l$ and we have
\begin{eqnarray*}
& &Y(\bar{u}^{(i)},x_{1})Y(\bar{u}^{(j)},x_{2})
=q_{ij}Y(\bar{u}^{(j)},x_{2})Y(\bar{u}^{(i)},x_{1}),\\
& &Y(\bar{v}^{(i)},x_{1})Y(\bar{v}^{(j)},x_{2})
=q_{ij}Y(\bar{v}^{(j)},x_{2})Y(\bar{v}^{(i)},x_{1}),\\
& &Y(\bar{u}^{(i)},x_{1})Y(\bar{v}^{(j)},x_{2})
-q_{ji}Y(\bar{v}^{(j)},x_{2})Y(\bar{u}^{(i)},x_{1})
=\delta_{ij}x_{1}^{-1}\delta\left(\frac{x_{2}}{x_{1}}\right)\ \ \ \
\ \ \ \
\end{eqnarray*}
for $i,j=1,\dots,l$. In view of Proposition \ref{puniversal-vq},
there exists a nonlocal-vertex-algebra homomorphism $\psi$ from
$V_{\bf Q}$ onto ${\rm Gr}_{F}(V(H,\S))$ such that
$$\psi (a^{(i)})=\bar{u}^{(i)},\ \ \ \psi(b^{(i)})=\bar{v}^{(i)}\ \
\ \mbox{ for }i=1,\dots,l.$$
  Since $V_{\bf Q}$ is simple and ${\rm Gr}_{F}(V(H,\S))\ne 0$,
  $\psi$ must be an
isomorphism. This proves ${\rm Gr}_{F}(V(H,\S))=V_{\bf Q}$.
 As $V_{\bf Q}$ is irreducible, by Proposition
\ref{p-filtration-irred}, $V(H,\S)$ is irreducible. Then $V(H,\S)$
is an irreducible quantum vertex algebra. As $V(H,\S)$ is a
nondegenerate quantum vertex algebra, it follows that $\S$ must be a
unitary rational quantum Yang-Baxter operator on $H$.
\end{proof}

\br{rold-theorem} {\em Here we correct an error in \cite{li-qva2}.
In Theorem 4.4 of \cite{li-qva2}, it was assumed that $H$ is a
finite-dimensional vector space equipped with a bilinear form
$\<\cdot,\cdot\>$ and $\S(x): H\otimes H\rightarrow H\otimes
H\otimes \C[[x]]$ is a linear map with $\S(0)=1$. Furthermore it was
assumed that ${\rm Gr}_{F}(V(H,\S))$ is linearly isomorphic to the
symmetric algebra over the space $H\otimes t^{-1}\C[t^{-1}]$ under a
canonical map. It was then claimed that $V(H,\S)$ is a nondegenerate
quantum vertex algebra.  The reasoning therein used Proposition 2.18
of \cite{li-qva2}, which claims that for a nonlocal vertex algebra
$V$ with a generating subset $T$, if ${\rm Gr}_{E}(V)$ is
nondegenerate, then $V$ is nondegenerate, where $E$ is the
increasing filtration of $V$, constructed in Lemma
\ref{lfiltration-E}. The discrepancy is that the filtration $F$ in
the assumption is different from $E$ in Proposition 2.18. Now, we
modify this theorem by adding the assumption that $\<\cdot,\cdot\>$
is a nondegenerate skew-symmetric bilinear form. This modified
theorem becomes a special case
 of Theorem \ref{tgeneral-zf} which uses Proposition
\ref{p-filtration-irred} of the current paper, instead of
Proposition 2.18 of \cite{li-qva2}. For the example that was
constructed in \cite{li-qva2}, $H=U\oplus U^{*}$ for some vector
space $U$, equipped with the standard nondegenerate skew-symmetric
bilinear form, so the modified theorem is applicable.}
 \er

For the rest of this section we consider the special case with
$\S(x)$ diagonal.

 \bd{ddaqx} {\em Let $l$
be a positive integer and let ${\bf Q}(x)=(q_{ij}(x))$ be an
$l\times l$ matrix over $\C[[x]]$ such that
\begin{eqnarray}\label{eqijqji}
q_{ij}(x)q_{ji}(-x)=1\;\;\;\mbox{ for }1\le i,j\le l.
\end{eqnarray}
Let $H$ be a $2l$-dimensional vector space with a basis
$\{a^{(1)},\dots, a^{(l)},b^{(1)},\dots,b^{(l)}\}$. Define a linear
map $\S(x): H\otimes H\rightarrow H\otimes H\otimes \C[[x]]$ by
\begin{eqnarray*}
& &\S(x) (a^{(i)}\otimes a^{(j)})=q_{ji}(x)(a^{(i)}\otimes
a^{(j)}),\ \ \S(x) (b^{(i)}\otimes b^{(j)})=q_{ji}(x)(b^{(i)}\otimes
b^{(j)}),\\
& &\S(x) (a^{(i)}\otimes b^{(j)})=q_{ij}(-x)(a^{(i)}\otimes
b^{(j)}),\ \ \S(x) (b^{(i)}\otimes
a^{(j)})=q_{ij}(-x)(b^{(i)}\otimes a^{(j)})
\end{eqnarray*}
for $i,j=1,\dots,l$, and equip $H$ with a bilinear form
$\<\cdot,\cdot\>$ defined by
\begin{eqnarray}
\<a^{(i)},a^{(j)}\>=\<b^{(i)},b^{(j)}\>=0\ \ \ \mbox{ and }\ \
\<a^{(i)},b^{(j)}\>=\delta_{ij}=-q_{ii}(0)\<b^{(j)},a^{(i)}\> .
\end{eqnarray}
 We define $V_{{\bf Q}(x)}$ to be the weak quantum
vertex algebra $V(H,\S)$ associated with this pair $(H,\S)$ from
${\bf Q}(x)$. That is, $V_{{\bf Q}(x)}$ is a universal weak quantum
vertex algebra with generators $a^{(i)},b^{(i)}$ $(i=1,\dots,l)$,
satisfying relations
\begin{eqnarray}\label{erelation3}
& &Y(a^{(i)},x_{1})Y(a^{(j)}(x_{2})
=q_{ij}(x_{2}-x_{1})Y(a^{(j)},x_{2})Y(a^{(i)},x_{1}),\nonumber\\
& &Y(b^{(i)},x_{1})Y(b^{(j)},x_{2})
=q_{ij}(x_{2}-x_{1})Y(b^{(j)},x_{2})Y(b^{(i)},x_{1}),\nonumber\\
& &Y(a^{(i)},x_{1})Y(b^{(j)},x_{2})
-q_{ji}(x_{1}-x_{2})Y(b^{(j)},x_{2})Y(a^{(i)},x_{1})
=\delta_{ij}x_{1}^{-1}\delta\left(\frac{x_{2}}{x_{1}}\right)\ \ \ \
\ \ \ \
\end{eqnarray}
for $i,j=1,\dots,l$. } \ed

With this notion, as an immediate consequence of Theorem
\ref{tgeneral-zf} we have:

 \bc{cnonconstant} Let $l$ be a positive
integer and let ${\bf Q}(x)=(q_{ij}(x))$ be an $l\times l$ matrix in
$\C[[x]]$ satisfying (\ref{eqijqji}). Suppose that $V$ is a nonzero
weak quantum vertex algebra with a generating subset
$T=\{u^{(i)},v^{(i)}\;|\; i=1,\dots,l\}$, satisfying the relations
\begin{eqnarray*}
& &Y(u^{(i)},x_{1})Y(u^{(j)}(x_{2})
=q_{ij}(x_{2}-x_{1})Y(u^{(j)},x_{2})Y(u^{(i)},x_{1}),\\
& &Y(v^{(i)},x_{1})Y(v^{(j)},x_{2})
=q_{ij}(x_{2}-x_{1})Y(v^{(j)},x_{2})Y(v^{(i)},x_{1}),\\
& &Y(u^{(i)},x_{1})Y(v^{(j)},x_{2})
-q_{ji}(x_{1}-x_{2})Y(v^{(j)},x_{2})Y(u^{(i)},x_{1})
=\delta_{ij}x_{1}^{-1}\delta\left(\frac{x_{2}}{x_{1}}\right)\ \ \ \
\ \ \ \
\end{eqnarray*}
for $i,j=1,\dots,l$. Then $V$ is an irreducible quantum vertex
algebra with $Gr_{F}(V)=V_{{\bf Q}(0)}$, where $F$ is the increasing
filtration of $V$, defined in Lemma \ref{linduction1}, associated
with the generating set $T$. Furthermore, such a nonzero quantum
vertex algebra $V$, if it exists, is unique up to isomorphism. \ec

The following is our second main result of this section:

\bt{tpre-diag} Let ${\bf Q}(x)=(q_{ij}(x))_{i,j=1}^{l}$ with
\begin{eqnarray}
q_{ij}(x)=q_{ij}p_{ij}(-x)/p_{ij}(x),
\end{eqnarray}
where $q_{ij}\in \C,\; p_{ij}(x)\in \C[[x]]$ satisfy that
$q_{ij}q_{ji}=1$, $p_{ij}(0)=1$  and $p_{ij}(x)=p_{ji}(x)$ for all
$i,j=1,\dots,l$. Then $V_{{\bf Q}(x)}$ is a (nonzero) irreducible
quantum vertex algebra with ${\rm Gr}_{F}(V_{{\bf Q}(x)})=V_{{\bf
Q}(0)}$.
 \et

\begin{proof} By Lemma \ref{lpseudo-auto},  there exist mutually commuting
pseudo-automorphisms $\Phi_{i}(x)$ for $1\le  i\le l$ of $V_{{\bf
Q}(0)}$ such that
\begin{eqnarray*}
& &\Phi_{i}(x)(u^{(j)})=u^{(j)}\otimes p_{ij}(x),\ \ \ \
\Phi_{i}(x)(v^{(j)})=v^{(j)}\otimes p_{ij}(x)^{-1} \ \ \ \mbox{ for
}j=1,\dots,l.\ \ \ \ \
\end{eqnarray*}
Set
\begin{eqnarray*}
a^{(i)}(x)=Y(u^{(i)},x)\Phi_{i}(x),\ \ \ \
b^{(i)}(x)=Y(v^{(i)},x)\Phi_{i}^{-1}(x)
\end{eqnarray*}
for $i=1,\dots,l$. We have
\begin{eqnarray*}
a^{(i)}(x_{1})a^{(j)}(x_{2})
&=&Y(u^{(i)},x_{1})\Phi_{i}(x_{1})Y(u^{(j)},x_{2})\Phi_{j}(x_{2})\\
&=&p_{ij}(x_{1}-x_{2})Y(u^{(i)},x_{1})Y(u^{(j)},x_{2})\Phi_{i}(x_{1})\Phi_{j}(x_{2})\\
&=&q_{ij}p_{ij}(x_{1}-x_{2})Y(u^{(j)},x_{2})Y(u^{(i)},x_{1})\Phi_{j}(x_{2})\Phi_{i}(x_{1})\\
&=&q_{ij}p_{ij}(x_{1}-x_{2})p_{ji}(x_{2}-x_{1})^{-1}a^{(j)}(x_{2})a^{(i)}(x_{1})\\
&=&q_{ij}(x_{2}-x_{1})a^{(j)}(x_{2})a^{(i)}(x_{1}),
\end{eqnarray*}
\begin{eqnarray*}
b^{(i)}(x_{1})b^{(j)}(x_{2})
&=&Y(v^{(i)},x_{1})\Phi_{i}^{-1}(x_{1})Y(v^{(j)},x_{2})\Phi_{j}^{-1}(x_{2})\\
&=&p_{ij}(x_{1}-x_{2})Y(v^{(i)},x_{1})Y(v^{(j)},x_{2})
\Phi_{i}^{-1}(x_{1})\Phi_{j}^{-1}(x_{2})\\
&=&q_{ij}p_{ij}(x_{1}-x_{2})Y(v^{(j)},x_{2})Y(v^{(i)},x_{1})
\Phi_{j}^{-1}(x_{2})\Phi_{i}^{-1}(x_{1})\\
&=&q_{ij}p_{ij}(x_{1}-x_{2})p_{ji}(x_{2}-x_{1})^{-1}b^{(j)}(x_{2})b^{(i)}(x_{1})\\
&=&q_{ij}(x_{2}-x_{1})b^{(j)}(x_{2})b^{(i)}(x_{1}),
\end{eqnarray*}
\begin{eqnarray*}
& &a^{(i)}(x_{1})b^{(j)}(x_{2})
-q_{ji}(x_{1}-x_{2})b^{(j)}(x_{2})a^{(i)}(x_{1})\\
&=&a^{(i)}(x_{1})b^{(j)}(x_{2})
-q_{ji}\frac{p_{ji}(x_{2}-x_{1})}{p_{ij}(x_{1}-x_{2})}b^{(j)}(x_{2})a^{(i)}(x_{1})\\
&=&p_{ij}(x_{1}-x_{2})^{-1}Y(u^{(i)},x_{1})Y(v^{(j)},x_{2})\Phi_{i}(x_{1})\Phi_{j}^{-1}(x_{2})
\\
& &-q_{ji}p_{ij}(x_{1}-x_{2})^{-1}Y(v^{(j)},x_{2})Y(u^{(i)},x_{1})
\Phi_{j}^{-1}(x_{2})\Phi_{i}(x_{1})\\
&=&\delta_{ij}x_{2}^{-1}\delta\left(\frac{x_{1}}{x_{2}}\right)
p_{ij}(x_{1}-x_{2})^{-1}\Phi_{j}^{-1}(x_{2})\Phi_{i}(x_{1})\\
&=&\delta_{ij}x_{2}^{-1}\delta\left(\frac{x_{1}}{x_{2}}\right)
p_{ij}(0)^{-1}\Phi_{j}^{-1}(x_{2})\Phi_{i}(x_{2})\\
&=&\delta_{ij}x_{2}^{-1}\delta\left(\frac{x_{1}}{x_{2}}\right).
\end{eqnarray*}
That is, we have defined an $(H,\S)$-module structure on $V_{{\bf
Q}(0)}$. Using the fact that $\Phi_{i}(x)({\bf 1})={\bf 1}$ for
$i=1,\dots,l$, we see that ${\bf 1}$ is a vacuum vector of the
$(H,\S)$-module $V_{{\bf Q}(0)}$. We claim that the $(H,\S)$-module
$V_{{\bf Q}(0)}$ is generated by ${\bf 1}$, so that $V_{{\bf Q}(0)}$
is a nonzero vacuum $(H,\S)$-module, which implies that $V_{{\bf
Q}(x)}\ne 0$. Let $E$ be the submodule of the $(H,\S)$-module
$V_{{\bf Q}(0)}$ generated by ${\bf 1}$. We have $\Phi_{i}(x){\bf
1}={\bf 1}$,
$$\Phi_{i}(x)a^{(j)}(x_{1})=p_{ij}(x-x_{1})a^{(j)}(x_{1})\Phi_{i}(x), \ \
\Phi_{i}(x)b^{(j)}(x_{1})=p_{ij}(x-x_{1})^{-1}b^{(j)}(x_{1})\Phi_{i}(x)$$
for $i,j=1,\dots,l$. It follows from induction that
$\Phi_{i}(x)E\subset E((x))$ for $1\le i\le l$. Similarly, we have
$\Phi_{i}(x)^{-1}E\subset E((x))$. As
$$Y(u^{(i)},x)=a^{(i)}(x)\Phi_{i}(x)^{-1},\ \ \
Y(v^{(i)},x)=b^{(i)}(x)\Phi_{i}(x),$$ it follows that $E$ is closed
under the components of $Y(u^{(i)},x)$ and $Y(v^{(i)},x)$ for
$i=1,\dots,l$. Consequently, $E=V_{{\bf Q}(0)}$, as claimed. This
completes the proof.
\end{proof}

\bex{example-l1} {\em Let $f(x)\in \C[[x]]$ with $f(0)=1$. In view
of Theorem \ref{tpre-diag}, there exists a (nonzero) irreducible
quantum vertex algebra $V$ with a generating set $\{a,b\}$ such that
\begin{eqnarray*}
& &Y(a,x_{1})Y(a,x_{2})=\pm
f(x_{1}-x_{2})f(x_{2}-x_{1})^{-1}Y(a,x_{2})Y(a,x_{1}),\\
& &Y(b,x_{1})Y(b,x_{2})=\pm f(x_{1}-x_{2})f(x_{2}-x_{1})^{-1}Y(b,x_{2})Y(b,x_{1}), \\
& &Y(a,x_{1})Y(b,x_{2})\mp
f(x_{2}-x_{1})f(x_{1}-x_{2})^{-1}Y(b,x_{2})Y(a,x_{1})
=x_{2}^{-1}\delta\left(\frac{x_{1}}{x_{2}}\right).
\end{eqnarray*}
Furthermore, any nonzero nonlocal vertex algebra $V$ with a
generating set $\{a,b\}$ satisfying the above relations is an
irreducible quantum vertex algebra and such a nonzero quantum vertex
algebra is unique up to isomorphism. } \eex

We furthermore consider certain nondegenerate, but not necessarily
irreducible, quantum vertex algebras.

\bp{phalf}  Let $q_{ij}\in \C,\;\;p_{ij}(x)\in \C[[x]]$ be such that
$q_{ij}q_{ji}=1$, $q_{ij}(0)=1$ and $q_{ij}(x)=q_{ji}(x)$ for $1\le
i,j\le l$. Suppose that $V$ is a weak quantum vertex algebra with a
generating subset $\{ u^{(1)},\dots,u^{(l)}\}$ such that
\begin{eqnarray}\label{erelation-uu}
Y(u^{(i)},x_{1})Y(u^{(j)},x_{2})
=q_{ij}p_{ij}(x_{1}-x_{2})p_{ij}(x_{2}-x_{1})^{-1}
Y(u^{(j)},x_{2})Y(u^{(i)},x_{1})
\end{eqnarray}
for $1\le i,j\le l$. Assume that $V$ has a normal basis in the sense
that the normal vectors $X_{1}\cdots X_{l}{\bf 1}$, where
$X_{i}=u^{(i)}_{-n_{1}}\cdots u^{(i)}_{-n_{r}}$ for $r\ge 0, \;
n_{i}\ge 1$ with $n_{1}\ge \cdots \ge n_{r}$ if $q_{ii}=1$ and
$n_{1}> \cdots > n_{r}$ if $q_{ii}=-1$, form a basis of $V$. Then
$V$ is a nondegenerate quantum vertex algebra. Furthermore, such a
quantum vertex algebra exists and is unique up to isomorphism. \ep

\begin{proof}
 Set ${\bf Q}(x)=(q_{ij}(x))_{i,j=1}^{l}$, where
 $q_{ij}(x)=q_{ij} p_{ij}(-x)p_{ij}(x)^{-1}$ for $i,j=1,\dots,l$.
By Theorem \ref{tpre-diag} and Corollary \ref{cnonconstant}, we have
a nondegenerate quantum vertex algebra $V_{{\bf Q}(x)}$ with ${\rm
Gr}_{F}(V_{{\bf Q}(x)})=V_{{\bf Q}(0)}$. Let $K$ be the subalgebgra
of $V_{{\bf Q}(x)}$ generated by $u^{(1)},\dots,u^{(l)}$. Then $K$
is nondegenerate and has a normal basis.

Now, we prove that $V$ is isomorphic to $K$. Using the pair $(H,\S)$
where $H=\coprod_{i=1}^{l}\C u^{(i)}$ equipped with the zero form
and $\S(x)(u^{(i)}\otimes u^{(j)})=q_{ij}(x)(u^{(i)}\otimes
u^{(j)})$ for $1\le i,j\le l$, we have a universal weak quantum
vertex algebra $E$ $(=V(H,\S))$ with generators
$u^{(1)},\dots,u^{(l)}$ and with defining relation
(\ref{erelation-uu}).  Using the commutation relation
(\ref{erelation-uu}) we see that $E$ is linearly spanned by all the
normal vectors. As $K$ is a homomorphic image of $E$, it follows
that $E$ also has a normal basis, so that $E$ is isomorphic to $K$.
The same reasoning shows that $E$ is isomorphic to $V$. It follows
that $E$, $K$ and $V$ are isomorphic and nondegenerate.
\end{proof}

\bex{yangian} {\em Let $V$ be a nonlocal vertex algebra with a
generating set $\{e,f,h\}$, satisfying the following relations:
 \begin{eqnarray*}
&&Y(e,x_{1})Y(e,x_{2})=-\frac{1+x_{1}-x_{2}}{1-x_{1}+x_{2}}Y(e,x_{2})Y(e,x_{1}),\\
&&Y(f,x_{1})Y(f,x_{2})=-\frac{1-x_{1}+x_{2}}{1+x_{1}-x_{2}}Y(f,x_{2})Y(f,x_{1}),\\
& &Y(h,x_{1})Y(h,x_{2})=Y(h,x_{2})Y(h,x_{1}),\\
& &Y(h,x_{1})Y(e,x_{2})
=-\frac{1+x_{1}-x_{2}}{1-x_{1}+x_{2}}Y(e,x_{2})Y(h,x_{1}),\\
& &Y(h,x_{1})Y(f,x_{2})
=-\frac{1-x_{1}+x_{2}}{1+x_{1}-x_{2}}Y(f,x_{2})Y(h,x_{1}),\\
 & &Y(e,x_{1})Y(f,x_{2})
=Y(f,x_{2})Y(e,x_{1}).
 \end{eqnarray*}
Assume that $V$ has a basis consisting of the vectors
\begin{eqnarray}
 e_{-m_{1}}\cdots e_{-m_{r}}h_{-n_{1}}\cdots
 h_{-n_{s}}f_{-k_{1}}\cdots f_{-k_{t}}{\bf 1}
 \end{eqnarray}
 for $r,s,t\ge 0,\; m_{1}>\cdots>m_{r}\ge 1, \ \
 n_{1}>\cdots >n_{s}\ge 1,\ \ k_{1}>\cdots >k_{t}\ge 1$.
By Proposition \ref{phalf}, $V$ is a nondegenerate quantum vertex
algebra.  This example is closely related to a quantum vertex
algebra associated to double Yangian $DY(sl_{2})$ in
\cite{li-infinity}.}
 \eex


\begin{thebibliography}{FKRW}

\bibitem [AB]{ab}
I. Anguelova and M. Bergvelt, $H_{D}$-Quantum vertex algebras and
bicharacters, arXiv:0706.1528.

\bibitem [BK]{bk}
B. Bakalov and V. Kac, Field algebras, {\em Internat. Math. Res.
Notices} {\bf 3} (2003) 123-159.

\bibitem [D1]{duncan1}
J. F. Duncan, Moonshine for Rudvalis's sporadic group I, arXiv:
math/0609449.

\bibitem [D2]{duncan2}
J. F. Duncan, Moonshine for Rudvalis's sporadic group II, arXiv:
math/0611355.

\bibitem[EK]{ek}
P. Etingof and D. Kazhdan, Quantization of Lie bialgebras, V, {\em
Selecta Math. (New Series)} {\bf 6} (2000) 105-130.

\bibitem[F]{fa}
L. D. Faddeev, Quantum completely integrable models in field theory,
{\em Soviet Scientific Reviews Sect. C} {\bf 1} (1980) 107-155.

\bibitem[FFR]{ffr}
A. Feingold, I. B. Frenkel and J. F. Ries, Spinor construction of
vertex operator algebras, triality, and $E_{8}^{(1)}$, Contemporary
Math. 121, 1991.

\bibitem[Ka]{kac}
V. Kac, {\it Vertex Algebras for Beginners}, University Lecture
Series {\bf 10}, Amer. Math. Soc., 1997.

\bibitem[Ku]{ku}
P. P. Kulish, Twist of quantum groups and noncommutative field
theory, arXiv: hep-th/0606056.

\bibitem[LL]{ll}
J. Lepowsky and H.-S. Li, {\em Introduction to Vertex Operator
Algebras and Their Representations}, Progress in Math.  {\bf 227},
Birkh\"auser, Boston, 2004.

\bibitem [Li1]{li-local}
H.-S. Li, Local systems of vertex operators, vertex superalgebras
and modules,  {\em J. Pure Appl. Alg.} {\bf 109} (1996) 143-195.

\bibitem [Li2]{li-twisted}
H.-S. Li, Local systems of twisted vertex operators, vertex
superalgebras and twisted modules, {\em Contemporary Math.} {\bf
193} (1996), 203-236.

\bibitem [Li3]{li-g1}
H.-S. Li, Axiomatic $G_{1}$-vertex algebras, {\em Commun. Contemp.
Math.} {\bf 5} (2003) 281-327.

\bibitem [Li4]{li-simple}
H.-S. Li, Simple vertex operator algebras are nondegenerate, {\em J.
Alg.} {\bf 267} (2003) 199-211.

\bibitem [Li5]{li-pseudo}
H.-S. Li, Pseudoderivations, pseudoautomorphisms and simple current
modules for vertex operator algebras, {\em Contemporary Math.} {\bf
392} (2005) 55-65.

\bibitem[Li6]{li-qva1}
H.-S. Li, Nonlocal vertex algebras generated by formal vertex
operators, {\em Selecta Math. (New Series)} {\bf 11} (2005) 349-397.

\bibitem[Li7]{li-gamma}
H.-S. Li, A new construction of vertex algebras quasi modules for
vertex algebras, {\em Advances in Math.} {\bf 202} (2006) 232-286.

\bibitem[Li8]{li-qva2}
H.-S. Li, Constructing quantum vertex algebras, {\em International
Journal of Mathematics} {\bf 17} (2006) 441-476.

\bibitem[Li9]{li-smash}
H.-S. Li, A smash product construction of nonlocal vertex algebras,
{\em Commun. Contemp. Math.} {\bf 9} (2007) 605-637.

\bibitem[Li10]{li-infinity}
H.-S. Li, Modules at infinity for quantum vertex algebras,
 submitted for publication; arXiv: 0705.0687.

\bibitem [Ma1]{manin1}
Yu. I. Manin, Quantum Groups and Non-commutative Geometry, CRM
Lecture Notes, 1989.

\bibitem [Ma2]{manin2}
Yu. I. Manin, Topics in Noncommutative Geometry, Princeton
University Press, 1991.

\bibitem [M]{mon}
S. Montgomery, Hopf Algebras and Their Actions on Rings, CBMS
Regional Conference Series in Mathematics, Number 82, AMS, 1993.

\bibitem [S]{s}
M. E. Sweedler, {\em Hopf Algebras,} Benjamin, New York, 1969.

\bibitem[MR]{mr}
M. Mintchev and E. Ragoucy, Interplay between Zamolodchikov-Faddeev
and reflection-transmission algebras, arXiv: math.QA/0306084.

\bibitem[W]{wei}
M. D. Weiner, Bosonic Construction of Vertex Operator Para-Algebras
from Sympletic Affine Kac-Moody Algebras, Ph.D thesis, State
University New York at Binghamton, 1994.

\bibitem[ZZ]{za}
A. B. Zamolodchikov and A. B. Zamolodchikov, Factorized $S$-matrices
in two-dimensions as the exact solutions of certain relativistic
quantum field theory models, {\em Ann. Phys.} {\bf 120} (1979)
253-291.

\end{thebibliography}
\end{document}